\documentclass[10pt, a4paper, twoside]{amsart}

\makeatletter
\def\@settitle{\begin{center}%
		\baselineskip14\p@\relax
		\normalfont\LARGE\scshape\bfseries
		\@title
	\end{center}%
}
\makeatother

\makeatletter

\def\subsection{\@startsection{subsection}{2}%
	\z@{.5\linespacing\@plus.7\linespacing}{.5\linespacing}%
	{\normalfont\large\bfseries}}

\def\subsubsection{\@startsection{subsubsection}{3}%
	\z@{.5\linespacing\@plus.7\linespacing}{.5\linespacing}%
	{\normalfont\itshape}}

\usepackage[usenames, dvipsnames]{color}
\definecolor{darkblue}{rgb}{0.0, 0.0, 0.45}

\usepackage[colorlinks	= true,
raiselinks	= true,
linkcolor	= darkblue, %MidnightBlue,
citecolor	= Mahogany,
urlcolor	= ForestGreen,
pdfauthor	= {Arman},
pdftitle	= {},
pdfkeywords	= {},
pdfsubject	= {},
plainpages	= false]{hyperref}

\usepackage{dsfont,amssymb,amsmath,subfigure, graphicx,enumitem,multirow} %,etaremune,savesym}
\usepackage{amsfonts,dsfont,mathtools, mathrsfs,amsthm,bm} 
\usepackage[amssymb, thickqspace]{SIunits}
\usepackage{algorithm,algorithmicx,fancyhdr}%
\usepackage{algpseudocode}
\allowdisplaybreaks
\date{\today}
\theoremstyle{theorem}
\newtheorem{Thm}{Theorem}[section]

\newtheorem{Lem}[Thm]{Lemma}

\newtheorem{As}[Thm]{Assumption}
\newtheorem{Prob}[Thm]{Problem}

\newtheorem{Def}[Thm]{Definition}

\newtheorem{Rem}[Thm]{Remark}

\theoremstyle{remark}

\newcommand{\ol}[1]{\overline{#1}}
\newcommand{\ul}[1]{\underline{#1}}

\newcommand{\set}[1]{\mathbb{#1}}

\graphicspath{{Fig/}}

\newcommand{\Acl}{A_{\text{cl}}}
\newcommand{\Set}[1]{\mathscr{#1}}

\newcommand{\Traj}[1]{\bm{\phi}^{\bm{u},\bm{#1}}(x)}

\newcommand{\Trajmpc}[1]{\bm{\phi}^{\bm{u}^{\text{mpc}},\bm{#1}}(x)}

\newcommand{\traj}[2]{\phi_{#2}^{\bm{u},\bm{#1}}(x)}

\newcommand{\trajmpc}[2]{\phi_{#2}^{\bm{u}^{\text{mpc}},\bm{#1}}(x)}
\newcommand{\trajmpC}[2]{\phi_{#2}^{\bm{u}^{\text{mpc}},\bm{#1}}}
\newcommand{\Dist}[3]{d_{#3}(#2,\mathcal{#1})}
\newcommand{\Seq}[2]{\{#1:#2\}}
\newcommand{\Targ}[1]{\mathbb{T}^{\mathbb{#1}}}
\newcommand{\targ}[2]{\mathcal{T}^{\mathbb{#1}}_{#2}}
\newcommand{\Targprod}[1]{\mathscr{T}^{\mathbb{#1}}_N}
\newcommand{\DisT}[4]{d_{#1}(#2,\targ{#3}{#4})}
\newcommand{\DisTM}[4]{d_{#1}(#2,\targ{#3}{#4}\ominus\tilde{L}_{#4}\hypreC{j})}

\newcommand{\Umpc}{\bm{u}^{\text{mpc}}(x)}
\newcommand{\UmpC}{\bm{u}^{\text{mpc}}}
\newcommand{\Umpci}[1]{u^{\text{mpc}}_{#1}(x)}
\newcommand{\UmpcI}[1]{u_{#1}^{\text{mpc}}}
\newcommand{\Umpcdyn}[2]{u^{\text{mpc}}_{#1}(#2)}
\newcommand{\ktrigx}{k^{\bm{w}}_{\text{trig}}(x)}

\newcommand{\Hyprec}{\bm{\mathcal{E}}(x)}
\newcommand{\HypreC}{\bm{\mathcal{E}}}
\newcommand{\hyprec}[1]{\mathcal{E}_{#1}(x)}
\newcommand{\hypreC}[1]{\mathcal{E}_{#1}}
\newcommand{\ermpc}[1]{e_{#1}^{\bm{w}}(x)}
\newcommand{\ermpC}[1]{e_{#1}^{\bm{w}}}
\newcommand{\Valopt}[1]{V^*_N(#1)}
\newcommand{\Valf}[2]{V_N(#1,#2)}

\newcommand{\Utilde}{\bm{\tilde{u}}(x;j)}

\newcommand{\Utildei}[1]{\tilde{u}_{#1}(x;j)}
\newcommand{\UtildeI}[1]{\tilde{u}_{#1}}
\newcommand{\Trajcan}{\bm{\phi}^{\bm{\tilde{u}},\bm{0}}(x;j)}

\newcommand{\Trajcani}[1]{\phi_{#1}^{\bm{\tilde{u}},\bm{0}}(x;j)}
\newcommand{\TrajcanI}[1]{\phi_{#1}^{\bm{\tilde{u}},\bm{0}}}
\newcommand{\Proj}[3]{\Pi_{#3}(#2,#1)}
\newcommand{\Slak}[1]{\bm{s}^{\mathbb{#1}}(x)}
\newcommand{\Slaki}[2]{s_{#2}^{\mathbb{#1}}(x)}
\newcommand{\Slaktilde}[1]{\bm{\tilde{s}}^{\mathbb{#1}}(x;j)}

\newcommand{\Slaktildei}[2]{\tilde{s}_{#2}^{\mathbb{#1}}(x;j)}
\newcommand{\SlaktildeI}[2]{\tilde{s}_{#2}^{\mathbb{#1}}}
\newcommand{\Ufeas}{\bm{u}^{\text{c}}(x;j)}
\newcommand{\UfeaS}{\bm{u}^{\text{c}}}

\newcommand{\UfeasI}[1]{u_{#1}^{\text{c}}}
\newcommand{\Xfeas}{\bm{\phi}^{\bm{u}^{\text{c}},\bm{0}}(x;j)}
\newcommand{\XfeaS}{\bm{\phi}^{\bm{u}^{\text{c}},\bm{0}}}

\newcommand{\XfeasI}[1]{\phi^{\bm{u}^{\text{c}},\bm{0}}_{#1}}
\newcommand{\Uhcan}{\bm{\hat{u}}^{\text{c}}(x;j+1)}
\newcommand{\UhcaN}{\bm{\hat{u}}^{\text{c}}}
\newcommand{\Uhcani}[1]{\hat{u}_{#1}^{\text{c}}}
\newcommand{\Xhcan}{\bm{\hat{\phi}}^{\bm{\hat{u}}^\text{c},\bm{0}}(x;j+1)}
\newcommand{\XhcaN}{\bm{\hat{\phi}}^{\bm{\hat{u}}^\text{c},\bm{0}}}

\newcommand{\XhcanI}[1]{\hat{\phi}_{#1}^{\bm{\hat{u}}^\text{c},\bm{0}}}
\title[A Decentralized Event-Based Approach for Robust Model Predictive Control]{
A Decentralized Event-Based Approach for Robust Model Predictive Control}

\author[A. Sharifi Kolarijani, S. C. Bregman , P. Mohajerin Esfahani, T. Keviczky]{Arman Sharifi Kolarijani, Sander C. Bregman, Peyman Mohajerin Esfahani, Tam\'{a}s Keviczky}
	\thanks{The authors are with the Delft Center for Systems and Control, TU Delft, The Netherlands ({\tt \{a.sharifikolarijani,s.c.bregman,p.mohajerinesfahani,t.keviczky\}@tudelft.nl}).}

\begin{document}
\maketitle

\begin{abstract}
In this paper, we propose an event-based sampling policy to implement a constraint-tightening, robust MPC method. The proposed policy enjoys a computationally tractable design and is applicable to perturbed, linear time-invariant systems with polytopic constraints. In particular, the triggering mechanism is suitable for plants with no centralized sensory node as the triggering mechanism can be evaluated locally at each individual sensor. 
From a geometrical viewpoint, the mechanism is a sequence of hyper-rectangles surrounding the optimal state trajectory such that robust recursive feasibility and robust stability are guaranteed. 
The design of the triggering mechanism is cast as a constrained parametric-in-set optimization problem with the volume of the set as the objective function. Re-parameterized in terms of the set vertices, we show that the problem admits a finite tractable convex program reformulation and a linear program relaxation. Several numerical examples are presented to demonstrate the effectiveness and limitations of the theoretical results.
\end{abstract}

%===============================================================================

\section{Introduction}
Nowadays, networked control systems (NCSs) generally demand an array of compatibility and efficiency measures from control design methods, such as utilization under shared resources, applicability to mobile tasks, and compatibility with digital communication infrastructures~\cite{baillieul2007control}. 
%The survey paper \cite{baillieul2007control} and the references therein provide an overview of these emerging challenges. 
Event-based control~(EBC) is a class of strategies that aim to improve efficiency of NCSs in the context of communication and computation. 
%Event-based approaches to execute control laws are a class of strategies that aim to systematically address the efficiency problem in the context of communication and computation. 
In EBC, the dynamics determine the instance to update a control action (contrary to the traditional case where a control action is updated \emph{periodically}) \cite{heemels2012introduction}.
%In EBC, the underlying dynamics determine the time at which it is required to update a control action (contrary to the traditional case in which the control action is updated in a \emph{periodic} fashion) \cite{heemels2012introduction}. 
There are two options to implement such an event-based logic: embedded in the sensory system, the so-called \emph{event-triggered control} \cite{tabuada2007event} and \cite{girard2015dynamic}, or embedded in the controller, the so-called \emph{self-triggered control} \cite{anta2010sample} and \cite{nowzari2012self}. The responsible entity to determine an update instance is known as the \emph{triggering mechanism}. 
%Following the wide-spread usage of event-based sampling policies in control methods, 
In particular, model predictive control (MPC) methods~\cite{camacho2013model} have been the subject of many studies in order to be amended with an EBC mindset. 
%In this regard, 
%and at least in a pure theoretical sense, 
%an area of control theory that has witnessed an increased level of interest 
%to employ event-based policies 
%is the class of model predictive control (MPC) methods \cite{camacho2013model}.  
%There are multiple reasons that encourage such a level of interest, ranging from inherent properties of MPC methods to practical advantages. 
%The computation of a control action is (often) a heavily involved task in MPC methods. One thus hopes by employing an event-based approach to reduce this burden.
%The process of computing the control action is (often) a heavily involved, computational  task in MPC methods. One thus hopes by employing an event-based implementation to reduce this computational burden. 
%In addition, the discarded predictions in a standard MPC setting can now serve as a basis to design a triggering mechanism.
%In addition, the predictions that are usually discarded in a standard MPC method can serve as a basis for the design of the triggering mechanism. 

% On the other hand, a natural drawback of MPC methods that hinders their utilization in practice is their required, extensive computational effort.
MPC methods are a class of \emph{on-line} optimization-based control approaches. 
In these methods, a measure of system performance is optimized over a finite horizon while states and inputs are subject to certain constraints. When the underlying dynamics is uncertain, the specific term \emph{robust} MPC (RMPC) is used for these methods in the literature~\cite{marruedo2002input}. 
We refer the interested reader to the survey papers \cite{mayne2000constrained} and \cite{mayne2014model} that discuss about different aspects of MPC. 

Traditionally, the controller solves the corresponding optimization problem at every time step and produces as outcomes two sequences of \emph{optimal} inputs and states. Then, the controller sends the first element of input sequence to the actuators and the remaining elements of the input sequence and the whole state sequence are discarded. These discarded predictions in a standard MPC setting can serve as a basis to design a triggering mechanism. 
%But on a positive note, as a by-product of control law computations, an educated guess about the possible future behavior of dynamics with its corresponding input trajectory is available. This extra information can be viewed as a basis to design a triggering mechanism. 
Moreover, the computational burden of MPC methods is a major drawback, hindering their usage in practice. 
%(by demanding more-advanced processing units or inability to handle plants with fast dynamics). 
One thus hopes by employing an EBC approach to reduce the frequency at which the underlying optimization problem is solved. 
Notice that there are already some techniques in the MPC literature (the so-called \emph{warm start} approaches~\cite{yildirim2002warm}) that exploit the computed sequences at the previous step to speed up the computation process. 

There is also a big incentive to exploit the computed sequences of MPC methods in a class of NCSs, namely, \emph{wireless sensor/actuator networks} (WSANs). In these systems, the most important concern is the energy efficiency, see e.g., \cite[Section IV.B]{willig2008recent}. 
%Roughly speaking, the energy depletion in a wireless node (e.g., placed on a sensor, actuator, or controller) is either from the micro-controller (responsible for logical/mathematical computations) or from the transceiver (responsible for sending and receiving data). 
The main source of energy depletion in a wireless node is the transceiver (responsible for sending and receiving data). 
To reduce the frequency of data transfer, it is hence more efficient (energy-wise) to aggregate the data into a single packet (if possible) and transmit the resulting packet at once over the communication network \cite{krishnamachari2002impact} and \cite{madden2002tag}. 
%Thus, it becomes clear that an event-based policy to implement an MPC approach is naturally compatible with this nature of WSANs. 

%On a practical note, there is also a big incentive to exploit the computed optimal trajectories. Particularly in \emph{wireless sensor/actuator networks} (WSANs), the single most important concern is the energy efficiency. Roughly speaking, the energy depletion in a wireless node (e.g., placed on a sensor, actuator, or controller) is either from the micro-controller (responsible for logical/mathematical computations) or from the transceiver (responsible for sending and receiving data). 
%The main source of energy depletion is however the transceiver. This observation is a key motivation to aggregate the data into a single packet (if possible) and transmit only the resulting packet over the communication network \cite{krishnamachari2002impact} and
%\cite{madden2002tag}. 
%Thus, it becomes clear that an event-based policy to implement an MPC approach is naturally compatible with this nature of WSANs. For a detailed discussion on this subject, we refer the reader to \cite{willig2008recent} and the references therein (see \cite[Subsection IV.B]{willig2008recent} for the energy efficiency topic).
%According to the above discussion, the combination of MPC methods with event-based implementation policies is in fact a (relevant) venue to explore.

\textbf{Statement of contribution:} In this paper, an event-triggered (ET) approach is proposed to implement an RMPC method on perturbed, linear time-invariant (LTI) systems.
%In this paper, an event-based sampling of an RMPC method is proposed for perturbed linear time-invariant (LTI) systems. f
The RMPC method is originally introduced in \cite{richards2006robust}. 
%The RMPC method is inspired by a constraint-tightening RMPC approach introduced in \cite{richards2006robust}. 
The core idea behind the ET approach is to construct a sequence of \emph{hyper-rectangles} around the optimal state sequence available from solving the RMPC problem. 
%(notice that these hyper-rectangles can be viewed as a certain type of \emph{weighted} $\ell_\infty$ norms, see Definition~\ref{def:hypR}). 
Then, these hyper-rectangles will be sent to the sensors. The optimal input sequence will also be transmitted to the actuators. 
Once the observed states at the sensory units leave these hyper-rectangles, a triggering happens and the states at the triggering instance will be transmitted to the controller. 
This procedure is then repeated in a sampled-data fashion. 
A key feature of the proposed ET approach is its ability to decide based on the local observation of each \emph{individual} sensor, whether to trigger or not. 
This feature stems from the fact that the sets describing the triggering mechanism are hyper-rectangles. Hence, the conditions required for a triggering in different states of the system are independent of each other. 
%To see this, note that this property is not guaranteed for general sets (even polytopic ones) describing a triggering mechanism. 
%the Euclidean coordinate axes are normal (or perpendicular) to the faces of the constructed hyper-rectangles. 
This feature is in particular appealing to systems with \emph{decentralized} (spatially dispersed) sensing units, including systems equipped with high-level (or supervisory) MPC methods, e.g., water treatment systems~\cite{sotomayor2002model}, HVAC systems~\cite{kelman2013analysis}, and commercial refrigeration systems~\cite{hovgaard2013nonconvex} to name a few. 
%One key feature of the event-based approach is its applicability to systems with decentralized sensing units. 
%(The term decentralized sensing refers to systems where the sensory units are spatially dispersed.) 
%The term \emph{decentralized sensing} refers to systems where the sensory units are spatially dispersed. 
%Thus, the ET approach can decide based on the local information, available at each individual sensor, whether to trigger or not.
Notice that the collocation of the triggering mechanism and the sensory units is physically impossible in such systems. Moreover, the addition of a \emph{central} node (on which the triggering mechanism is placed on) to collect the sensory data comes at the price of extra communication bandwidth usage.  
%Although the sensory units' information can be collected on a single \emph{centralized} node to place a triggering mechanism in such systems, we make the general assumptions that (i) the extra burden on the communication component is undesirable, and (ii) the complications that arise from the non-collocation of the triggering mechanism with the sensory units reduce the practicality of the triggering mechanism. 
On the theoretical side, the design of the ET approach is \emph{decoupled} from the design of the underlying RMPC method. 
As a result, a fair comparison between the performances of the ET and standard implementations of the RMPC method becomes possible. 
%Another important feature of the proposed ET approach is the fact that the triggering mechanism's design is \emph{decoupled} from the design of the underlying RMPC method. 
This paper extends the results of the authors' previous work in \cite{bregman2017robust} in multiple directions, in particular, by simplifying the triggering ``{law}". The approach in \cite{bregman2017robust} requires an ``{advanced}" triggering mechanism that is responsible for (i)~constructing certain input and state sequences, (ii) evaluating the satisfaction of MPC's constraints by these sequences, and (iii) comparing the values of the cost function based on the constructed sequences with the value function at the last triggering instance. 
%This paper extends the results of the authors' previous work in \cite{bregman2017robust} in multiple directions. A set-regulator RMPC problem with aperiodic actuation is proposed in \cite{bregman2017robust} with a collocated controller and sensory units setup. The main limitation of the approach proposed in \cite{bregman2017robust} is the requirement of an ``{intelligent}" sensory system.}  
%We borrow the formulation of the RMPC method from \cite{bregman2017robust}.
The main contributions of the paper are summarized as follows.
\begin{itemize}
\item {\bf Decoupled recursive feasibility and stability:} Given an RMPC method in place, we propose a set-theory-based, ET approach that preserves robust recursive feasibility and robust stability. 
%Unlike the existing literature, 
The proposed approach is decoupled from the control synthesis process and does not require additional assumptions, such as extra conditions on eigenvalues of weighting matrices in the cost function or the need to define user-specified thresholds for the triggering mechanism (Theorem~\ref{theo_4}).
\item {\bf Decentralized applicability:} The proposed approach enjoys a decentralized triggering mechanism that only requires local sensory information (Definition~\ref{def:tm}).
\item \textbf{Tractable convex program reformulation:} We show that a certain type of non-convex volume-maximization problem with set-based constraints that is deployed to design the triggering mechanism admits a finite tractable convex program (CP) reformulation (Theorem~\ref{theo:lp-nlp-Box_1}). 
%Finally, we introduce a heuristic CP framework to address the limitations of the proposed CP and LP frameworks (Remark~\ref{rem:reduc_dir}).
\item \textbf{Suboptimal linear program relaxation:} Motivated by an approach in the literature, we further show that a linear program (LP) relaxation of the CP reformulation is possible (Theorem~\ref{theo:lp-nlp-Box_2}).
\end{itemize}

\textbf{Literature review:} In what follows, we first review several event-triggered, MPC approaches. 
We then close this section by giving a brief account of several computationally efficient approaches that are customized for MPC problems.
%In what follows, we first review several approaches that are closely related to the problem considered in this paper. 
%We then give a brief account of algorithmic-oriented approaches in the literature that reduce the computational complexity of their MPC methods with customized algorithms.
\newline
\emph{Related works:} 
%In order to avoid repetition of terms (unless the contrary is mentioned otherwise), 
Let us first mention the shared properties of the references below: linear discrete-time models, event-triggering mechanisms, constrained MPC methods, minimal (to none) coupling of the parameters of the triggering mechanism and the considered MPC method, and a computationally viable approach to design the triggering mechanism. 
%{\color{blue}Let us now elaborate the reasoning behind the last two properties. First, in an \emph{ideal} case, one seeks the possibility of a complete \emph{decoupling} between the parameters of the triggering mechanism and of the considered MPC method. (By doing so, a fair comparison between the performances of the event-based and standard implementations of an MPC method will be possible.) After all, our main goal is to provide an implementation policy for an MPC method with the awareness of communication issues. Hence, we only mention the approaches that have minimal to no inter-connections between the two sets of design parameters. Second, we also do not consider event-triggering approaches that are computationally more expensive compared to the underlying MPC problem (e.g., approaches that require to solve some type of an integer program). As mentioned above, the prohibitive computational requirement of MPC methods is the main factor that limits their application in practice. Thus, the viability of such complicated event-triggering designs becomes questionable in practice. The review of related works follows.} 

%Goodwin[2004]:
To deal with practical issues such as a band-limited communication channel, a novel design approach for NCSs is proposed in
 \cite{goodwin2004moving}. They employ the notion of \emph{moving horizon} \cite{rao2003constrained} to design the estimator and controller. A remarkable character of their approach is its ability to decide \emph{on-the-fly} which input channel should be updated (i.e., a certain type \emph{input-channel} event-triggering control).
%[2010]:
In case of collocated controller and actuator units, an event-based estimator with a bounded covariance matrix is designed in \cite{sijs2010integration}.
While the estimator receives data via a Lebesgue sampling approach, it periodically updates the controller's information regarding the disturbances with a polytopic over-approximation of the covariance matrix.
%[2012]:
The authors of \cite{bernardini2012energy} propose an interesting transmission strategy for wireless sensor/controller communications with practical energy-aware provisions 
(the controller is collocated with the actuator system). Using some predefined thresholds for each state's sensor (i.e., an $\ell_1$-type triggering mechanism), the controller is computed offline using an \emph{explicit} MPC approach \cite{bemporad2002explicit}.  
%[2013]:
Based on a prescribed 2-norm ball around the optimal state trajectory, the authors in \cite{lehmann2013event} propose a triggering mechanism for WSANs. They show that the approach is robustly stable to a set that is a function of the radius of threshold ball and the maximal 2-norm of disturbance. 
%[2014]:
For linear, continuous-time dynamical systems affected by a Wiener process, a co-design method (i.e., simultaneous design of the scheduler and the controller) is proposed in \cite{antunes2014rollout}. The main idea is inspired by the notion of \emph{rollout} from \emph{dynamic programming} \cite{bertsekas2005dynamic}. More importantly, the authors show that under some mild conditions, an event-based control approach outperforms a traditional control approach w.r.t. closed-loop performance/average transmission rate. (Notice that for most of the approaches in the literature including our paper such a guarantee is not provided.)  
%[2017]:
A set theoretic triggering mechanism is introduced in \cite{brunner2017robust} for systems with collocated controller and sensory units. The approach is inspired by the \emph{tube-based} MPC proposed in \cite{rakovic2012homothetic}. By exploiting the known probability distribution of disturbance, they also guarantee an average sampling rate. However, their tube-contraction method requires a certain type of realization of a discrete-time system, see \cite[Remark 8]{brunner2017robust}.
%[2017]: 
Demirel et al., introduce a sensor/actuator event-triggering mechanism for control systems with limited number of control messages (i.e., communication and computation resources are scarce) \cite{demirel2017optimal}. They relax the underlying combinatorial problem into a convex one by an appropriate definition of event thresholds. 
In \cite{incremona2017asynchronous}, a packetized approach is proposed for input-affine, nonlinear systems with bounded additive disturbances in continuous-time. In the proposed approach, an RMPC controller (connected via a communication network to the plant) takes into account the mismatched uncertainties while an integral sliding-mode controller~\cite{utkin1996integral} (placed at the plant) counters the effect of the matched uncertainties.  
\newline
\emph{Algorithmic viewpoint:} 
%The reasoning behind this algorithmic viewpoint is as follows. 
An MPC optimization problem is computationally expensive by itself. 
%(the evidence is the substantial body of work that has been done to customize algorithms to MPC problems). 
Hence, the merit of an event-based policy of implementation would be lost if the mechanism demands a drastically higher computational effort compared to the underlying MPC problem. 
%We should highlight that all the approaches discussed next are related to linear time-invariant (LTI) systems.
Dunn and Bertsekas in \cite{dunn1989efficient} exploit the structure of their problem to reduce the cubic complexity of computing a Newton step to a linear one. 
In \cite{wang2010fast}, the authors use a specific \emph{ordering} of decision variables to promote a sparse structure that decreases the cost of computing a control action. 
The authors in \cite{richter2012computational} employ a simple, gradient-based algorithm to solve an MPC  problem while providing \emph{a priori} computational complexity certificate.

The layout of the paper is as follows. The mathematical notions used in the paper are outlined in Section~\ref{sec:pre}. Section~\ref{sec:st_RMPC} is devoted to the considered RMPC method. The main results regarding the event-based implementation policy are introduced in Section~\ref{sec:eb_RMPC}. Section~\ref{sec:proof} contains the technical proofs. Several numerical examples are presented in Section~\ref{sec:exp} to evaluate the effectiveness and limitations of the theoretical results. Finally, we present several future research directions in Section~\ref{sec:conc}.

\section{Notation and Preliminaries}
\label{sec:pre}

We begin with a brief review of the mathematical preliminaries employed in the rest of the paper.

\textbf{Notation:} The set of non-negative integers is denoted by $\mathbb{Z}_{\geq 0}$. 
Given positive integers $m$ and $n$, $\mathbb{R}^m$ and $\mathbb{R}^{m\times n}$ represent the $m$-dimensional Euclidean space and the space of $m\times n$ matrices with real entries, respectively. 
Given two integers $i,j$ where $i\leq j$, $\Seq{i}{j}:=\{i,i+1,\ldots,j\}$. 
%Given a vector $v\in\mathbb{R}^n$, $v^i$ represents the $i$-th entry of $v$. 
For any pairs of vectors $a,b\in\mathbb{R}^n$, the inequality $a<(\leq)b$ is realized in a component-wise manner. 
Given a vector $v\in\mathbb{R}^n$ and a scalar $p\geq1$, $\|v\|_p$ denotes the $p$-norm $\big(\sum_{i=1}^n~(v^i)^p  \big)^{1/p}$. 
%, i.e., $a^i<(\leq)b^i$, for all $i\in\Seq{1}{n}$. 
Given a matrix $M\in \mathbb{R}^{m \times n}$, $M_{ij}$ denotes the $i$-th row, $j$-th column entry of $M$. Moreover, the matrix $M^+\in\mathbb{R}^{m\times n}$ is the matrix with entries~$M^+_{ij}:=\max\{0,M_{ij} \}$. 
%A positive definite matrix $M$ is denoted by $M\succ 0$. 
The $n\times n$ zero and identity matrices are denoted by $\mathsf{0}_n$ and $\mathsf{I}_n$, respectively. 
Given a set $\mathcal{S}\subset \mathbb{R}^n$ and a matrix $M\in\mathbb{R}^{m\times n}$, the set $M\mathcal{S}$ denotes the set $\{c\in\mathbb{R}^m:~\exists s\in\mathcal{S}, Ms=c\}$. 
%Given two sets $\mathcal{A}$ and $\mathcal{B}$ in $\mathbb{R}^n$, $\mathcal{A}/\mathcal{B}:=\{x\in\mathcal{A}:~x\notin\mathcal{B}\}$. 
Given a matrix $M \succ 0$ (i.e., positive definite), the squared weighted distance of a point $r\in\mathbb{R}^n$ from a closed set $\mathcal{S}\subset \mathbb{R}^n$ is defined as
%	\begin{align*}
%	\label{MPC:distance_metric}
$	\Dist{S}{r}{M} := \min_{s \in \mathcal{S}}~ ||r-s||^2_M = \min_{s \in \mathcal{S}}~ (r-s)^\top M (r-s)$. 
Denote the projection of $r$ onto $\mathcal{S}$ by $\Proj{\mathcal{S}}{r}{M} \in\text{argmin}_{s\in\mathcal{S}} \Dist{S}{r}{M}$. Note that when $\mathcal{S}$ is also convex, the projection is unique. 
Given sets $\mathcal{C}$ and $\mathcal{D}$, the Pontryagin difference $\mathcal{C}\ominus \mathcal{D}$ and the Minkowski sum $\mathcal{C}\oplus\mathcal{D}$ are defined as
%\begin{align*}
%\label{eq:pontDif}
$\mathcal{C}\ominus \mathcal{D}  :=\{c:~c+d\in \mathcal{C}, \forall d\in \mathcal{D} \}$ and $
\mathcal{C}\oplus \mathcal{D}  :=\{c+d:\forall c\in \mathcal{C},\forall d\in \mathcal{D} \}$, respectively.  
The function $\text{sign}(\cdot)$ represents the standard sign function. 
%For a set~$\mathcal{S}\subset\mathbb{R}^n$, the \emph{characteristic} function~$\Ind{\mathcal{S}}{s}=1$ if $s\in\mathcal{S}$ and $\Ind{\mathcal{S}}{s}=0$ otherwise. 
Given a set~$\mathcal{X}\in\mathbb{R}^n$ and an \emph{extended} real-valued function~$f:\mathcal{X}\to[-\infty, +\infty]$, the \emph{effective domain} of $f$ is the set $\text{dom}(f)=\{x\in\mathcal{X}:f(x)<\infty\}$. 
%Moreover, $f$ is called \emph{proper} if $f(x)< \infty$ for some $x\in\mathcal{X}$ and $f(x)>-\infty$ for all $x\in\mathcal{X}$.

%The section is divided to two parts: RMPC related notions and Set theory related notions.
%\end{align*}
%\end{Def}
The following result will be used frequently in the development of the triggering mechanism.
\begin{Lem}[Set-difference lower bound \cite{richards2006robust}]
\label{lem:lem_1}
Let $r$ be a vector in $\mathbb{R}^n$, $\mathcal{B}$ and $\mathcal{C}$ be two compact sets in $\mathbb{R}^n$, and $M$ be a positive definite matrix in $\mathbb{R}^{n \times n}$. Then, 
%\begin{align*}
%\label{MPC:distance_metric_upperbound}
$\Dist{B}{r+c}{M} \leq d_{M}(r, \mathcal{B}\ominus\mathcal{C})$, for all $c \in \mathcal{C}$.
%\end{align*}
\end{Lem}

We now revisit some notions from convex analysis (see e.g., \cite[Section~2]{kolmanovsky1998theory} for a compact exposition of the subject). 
%\begin{Def}[Support function]
%\label{def:supF}
Given a set $\mathcal{S} \subset\mathbb{R}^n$, the support function of $\mathcal{S}$ evaluated at $\eta\in\mathbb{R}^n$ is 
%\begin{align*}
%\label{eq:supF}
$h_\mathcal{S}(\eta)  :=  \sup_{s\in \mathcal{S}} ~ \langle \eta, s\rangle$. 
%\end{align*}
%\end{Def}
The domain $\mathcal{K}_{\mathcal{S}}$ on which the support function is defined is a convex cone pointed at the origin. 
If $\mathcal{S}$ is bounded, then $\mathcal{K}_\mathcal{S}:=\mathbb{R}^n$. 
Given a matrix $M\in\mathbb{R}^{n\times m}$ and a vector $v\in \mathbb{R}^n$, if $M^{\top} v \in \mathcal{K}_\mathcal{S}$, then 
%\begin{align*}
%\label{eq:supM}
$h_{M\mathcal{S}}(v):=h_\mathcal{S}(M^\top v)$.
%\end{align*}
Suppose $\mathcal{S}\subset \mathbb{R}^n$ is closed and convex. Then, 
%\begin{align*}
%\label{eq:hypPlane}
$\mathcal{S}:=\{s\in\mathbb{R}^n: \langle \eta ,s\rangle \leq h_\mathcal{S}(\eta),\forall \eta \in \mathcal{K}_{\mathcal{S}} \}$, 
%\end{align*}
i.e., the intersection of its supporting halfplanes. 
%\begin{Def}[Polyhedron]
%\label{def:hPoly}
A set $\mathcal{S}\subset \mathbb{R}^n$ is called a polyhedron, if 
%\begin{align*}
%\label{eq:polyH}
$\mathcal{S}= \{s\in \mathbb{R}^n: A_\mathcal{S} s \leq b_\mathcal{S} \},$ $A_\mathcal{S}\in\mathbb{R}^{m\times n}$, $b_\mathcal{S}\in\mathbb{R}^m$.
%\end{align*}
%\end{Def}
If the polyhedron $\mathcal{S}$ is bounded, the set is called a \emph{polytope} and its representation given above is known as the \emph{H-representation}. Furthermore, the support function $h_\mathcal{S}(\eta)$ of a polytope~$\mathcal{S}$ is the solution of the LP, 
%\begin{equation*}
%\label{eq:subFP}
%\begin{aligned}
$h_\mathcal{S}(\eta) =  \max_s ~ \langle \eta,s \rangle$ subject to  $A_\mathcal{S} s \leq b_\mathcal{S}$. 
%\end{aligned}
%\end{equation*}
Given the \emph{H}-representation of a polytope, we employ the notations $a_{i,\mathcal{S}}\in\mathbb{R}^{1\times n}$ and $a_{\mathcal{S},j}\in\mathbb{R}^{m\times 1}$ to denote the $i$-th row and the $j$-th column of $A_{\mathcal{S}}$, respectively. Moreover,  $b_{i,\mathcal{S}}$ is the $i$-th entry of $b_{\mathcal{S}}$. 
Given a polyhedron $\mathcal{S}\subset \mathbb{R}^n$ and a set $\mathcal{V}\subset \mathbb{R}^n$, 
assume that $h_{\mathcal{V}}(a_{i,\mathcal{S}}^\top)$ is well-defined for all $i\in \Seq{1}{m}$. Then, 
%\begin{multline}
%\label{eq:pontDH}
$\mathcal{S}\ominus \mathcal{V}:= \big\{z\in\mathbb{R}^n:~ 
 \langle a_{i,\mathcal{S}}^\top,z  \rangle \leq b_{i,\mathcal{S}} - h_{\mathcal{V}}(a_{i,\mathcal{S}}^\top), \forall i\in \Seq{1}{m} \big\}$. 
%\end{multline}
%\begin{Def}[Full-dimensional hyper-rectangle]
%\label{def:hypR}
For any vector-pairs $l,u\in\mathbb{R}^n$ such that $l<u$, the full-dimensional convex polytope
%\begin{align*}
$\mathcal{B}(l,u):=\{x\in\mathbb{R}^n:~l\leq x \leq u  \}
 =\{x\in\mathbb{R}^n:~ A_{\mathcal{B}}x \leq b_{\mathcal{B}}  \}$
%\end{align*}
is called a \emph{hyper-rectangle}, where $A_{\mathcal{B}}:=[\mathsf{I}_n~-\mathsf{I}_n]^\top$ and $b_{\mathcal{B}}=[u^\top~-l^\top]^\top$.
%\end{Def}

\section{Robust model Predictive Control Method}
\label{sec:st_RMPC}

In this section, we introduce the class of constrained dynamical systems considered in this paper, followed by the description of the RMPC method. 
At last, we formally state the problem addressed in this paper.

Consider an LTI system with a bounded additive disturbance given by
\begin{align}
\label{eq:dynamics}
x^+= A x + B u + w,
%x_{k+1} = A x_k + B u_k + w_k,~\text{for all}~ k\in\mathbb{Z}_{\geq0},
\end{align}
where $x^+$ is the successor state and $x$, $u$, and $w$ are the current state, input and disturbance, respectively. The current state, input, and disturbance are subject to the hard constraints
\begin{align}
\label{eq:constraints}
x \in \mathbb{X}\subset \mathbb{R}^{n_x},~
u \in \mathbb{U}\subset \mathbb{R}^{n_u},~
w \in \set{W} \subset \mathbb{R}^{n_x}.
\end{align}
A system is called the \emph{nominal} system associated with \eqref{eq:dynamics} when $w=0$. 
%\begin{align}
%\label{eq:nom_dynamics}
%x^+ = A x + B u.
%\bar{x}_{k+1} = A \bar{x}_k + B u_k,~\text{for all}~ k\in\mathbb{Z}_{\geq0}.
%\end{align}
Given a positive integer $N$, let $\Set{U}:=\mathbb{U}^N=\prod_{i=0}^{N-1}\mathbb{U}$ ($\mathscr{W}:=\set{W}^N$) denote the class of admissible control sequences $\bm{u}:=\{u_i\}_{i\in\Seq{0}{N-1}}$ (admissible disturbance sequences $\bm{w}:=\{w_i\}_{i\in\Seq{0}{N-1}}$). 
Initiated at state~$x$, the solution to \eqref{eq:dynamics} at time $i$ with the control and disturbance sequences~$\bm{u}$ and $\bm{w}$, respectively, is denoted by $\traj{w}{i}$. 
Similarly, we define $\Traj{w}:=\{\traj{w}{i}\}_{i\in\Seq{0}{N}}$. 
Moreover, let $\traj{0}{i}$ denote the nominal solution with the input sequence $\bm{u}$ initiated at state~$x$. The RMPC method is designed such that the state~$x$ and the input~$u$ eventually converge to some user-defined \emph{target sets} $\Targ{X} \subset \mathbb{R}^{n_x}$ and $\Targ{U} \subset \mathbb{R}^{n_u}$, respectively, while the constraints~\eqref{eq:constraints} are satisfied at all times.  

\begin{As}[System \& constraint sets]
\label{as:ass_0}
(i)~Nominal controllability: The pair $(A, B)$ is controllable.
(ii)~Polytopic sets: The sets $\mathbb{X}$, $\mathbb{U}$, $\Targ{X}$, $\Targ{U}$, and $\set{W}$ are all convex, compact polytopes containing their underlying spaces' origin in their interior.
\end{As}

We start with introducing two types of feedback gains which are used in the RMPC method and are essential for the construction of the triggering mechanism. 
Let $F\in\mathbb{R}^{n\times m}$ be a given feedback gain that guarantees the stability of the nominal system with $u=Fx$. 
The \emph{nominal} gain~$F$ can be designed so that a satisfactory performance (e.g., in an LQ optimal control sense) is guaranteed for the nominal system.

Let integer $N\geq n_x +1$ be the horizon length of the RMPC method and integer $M$ be given, where $ M\in\Seq{n_x}{N-1}$. 
Suppose next that a set of feedback gains $\bm{K}=\{K_i\}_{i\in\Seq{0}{N-1}}$ are given such that $\prod_{i=1}^M(A+BK_i)=0$, i.e., for all $k\geq M$, $\traj{0}{k}=0$. 
%the input sequence $\mathbf{u}=\{K_i x^{\mathbf{u}, \mathbf{w}}(i;x)\}_{i=0}^{N-1}$ renders the nominal dynamics~\eqref{eq:nom_dynamics} nilpotent in $M$-steps
We call the set of gains~$\bm{K}$ the \emph{tightening gains} since these gains are employed in the state and input constraint tightening process. 
We refer the interested reader to \cite[Section~IV]{richards2006robust} for a possible approach to construct the gains~$\bm{K}$. 
The constraint tightening approach is applied to the input, state, input target, and state target sets, that is, for all $ i \in \Seq{0}{N-2}$,
\begin{subequations}\label{tightening}
	\begin{alignat}{2}
	\label{tight:input} \mathcal{U}_0 &= \mathbb{U}, \quad &\mathcal{U}_{i+1} &= \mathcal{U}_i \ominus K_i L_i \set{W},\\
	\label{tight:state} \mathcal{X}_0 &= \mathbb{X}, \quad &\mathcal{X}_{i+1} &= \mathcal{X}_i \ominus L_i \set{W},\\
	 \targ{U}{0} &= \Targ{U}, \quad & \targ{U}{i+1} &=  \targ{U}{i} \ominus K_i L_i \set{W},\\
	 \targ{X}{0} &= \Targ{X}, \quad & \targ{X}{i+1} &=  \targ{X}{i} \ominus L_i \set{W},
	\end{alignat}
\end{subequations}
where $L_0=\mathsf{I}_{n_x}$ and $L_{i+1} = (A + BK_i)L_i$ for all $i\in\Seq{0}{N-2}$. 
Notice that the $M$-step nilpotency of the set of gains~$\mathbf{K}$ implies that for all $i\in\Seq{M}{N-1}$, $L_i=\mathsf{0}_{n_x}$.

Let the terminal set~$\mathcal{X}_f\subset \mathbb{R}^{n_x}$ be a \emph{control invariant set} for the nominal system, i.e., $(A + BF) \xi \in \mathcal{X}_f$ for all $\xi\in\mathcal{X}_f$. 
\begin{As}[Terminal set]
\label{as:terminal}
For all $\zeta \in \mathcal{X}_f$, the following conditions hold:
	\begin{align*}
  \zeta \in \mathcal{X}_{N-1} \cap \targ{X}{N-1},~ 
 F \zeta \in \mathcal{U}_{N-1} \cap \targ{U}{N-1}.
	\end{align*}
\end{As}
For the sake of notational simplicity, let us define $\mathscr{U}_N:=\prod_{i=0}^{N-1}\mathcal{U}_i$ and $\mathscr{X}_N:=\prod_{i=0}^{N-1}\mathcal{X}_i\times\mathcal{X}_f$.
%\textbf{(3) Cost function:} 
%Let us denote the input trajectory $\{u_{k+i|k}\}_{i=0}^{N-1}$ and the state trajectory $\{x_{k+i|k}\}_{i=0}^{N}$ by $\mathbf{U}_{k|k}$ and $\mathbf{X}_{k|k}$, respectively. 
The \emph{cost function} of the RMPC problem is
\begin{align}
\label{eq:cost_mpc}
V_N(x, \bm{u}) := \sum_{i = 0}^{N-1}  d_{Q}(\traj{0}{i},\targ{X}{i} ) + d_{R}(u_i, \targ{U}{i})
 +\delta_{\text{feas}}\big(\bm{u}, \Traj{0}\big),
\end{align}
%We are now in a position to introduce the RMPC problem.
where $\delta_{\text{feas}}\big(\bm{u}, \Traj{0}\big)=0$ if $\bm{u} \in \Set{U}_N$ and $\Traj{0}\in\Set{X}_N$, and $=\infty$ otherwise, is the \emph{indicator} function of the set~$\Set{U}_N\times \Set{X}_N$. Notice that the input and state constraints are embedded in the objective function via the indicator function. The optimization problem for a finite horizon~$N$ with an initial state~$x$ reads as
\begin{align}
\label{MPC_tightened_target}
V_N^*(x) := \min_{\bm{u}} V_N(x, \bm{u}),
\end{align}
with $\Umpc := \text{argmin}_{\bm{u}} V_N(x, \bm{u})$ as the optimal input sequence. 
When it is clear from the context, we may instead use the shorthand notation~$\UmpC$. 
The above sequence of inputs is indeed an optimal solution to a nominal (i.e., $\bm{w} = \bm{0}$) finite optimization problem emerging in the context of finite horizon MPC in the rest of the paper. In this light, we denote this nominally optimal controller by a similar label, for which the associated nominal state sequence is $\Trajmpc{0}$. 
%Let $\text{dom}(V^*_N)$ be the effective domain of $V_N^*(\cdot)$.} 
%(with the standard convention that the value function~$V_N^*(x):=+\infty$ if the problem~\eqref{MPC_tightened_target} is infeasible).} 
%The domain~$\text{dom}(V^*_N)$ of the value function is
%\begin{align}
%\label{eq:val_dom}
%\big\{x\in\mathbb{R}^{n_x}: \Traj{0}\in\Set{X}_N \text{~with~}\bm{u} \in \mathscr{U}_N, \big\}.
%\end{align}

%\begin{Rem}[Standard implementation]
In a standard RMPC setting, the optimal control problem~\eqref{MPC_tightened_target} is solved. The first element $\Umpci{0}$ of $\Umpc$ is then applied to the plant yielding to the closed-loop dynamics $
x^+=Ax+B\Umpci{0}+w$. 
In an event-based setting, the triggering mechanism generally exploits the optimal state sequence~$\Trajmpc{0}$ in order to possibly employ the rest of elements in the nominally optimal input vector $\Umpc$. 
The challenge in designing the triggering mechanism is then to guarantee robust stability and robust recursive feasibility of the resulting event-triggered, closed-loop dynamics. 
%\end{Rem}

\begin{Def}[Triggering mechanism]
\label{def:tm}
Given an initial state~$x$ and a sequence of (possibly) state-dependent, hyper-rectangular sets $\Hyprec:=\hypreC{0}\cup \{\hyprec{i}\}_{i=1}^{N-1}\subset (\mathbb{R}^{n_x})^{N}$, the triggering instance is defined by 
\begin{align}
\label{eq:trig_state}
\ktrigx:=\min\big\{j\in\Seq{0}{N-1}:
\trajmpc{w}{j} - \trajmpc{0}{j} \notin \hyprec{j} \big\},
\end{align}
where $\hypreC{0}:=\mathbb{R}^{n_x}$.
\end{Def}
The quantity~$\ktrigx$ is known as the \emph{inter-execution} time in the literature. One can observe that $\trajmpc{w}{0}= \trajmpc{0}{0}=x$. As a result, $\trajmpc{w}{0}- \trajmpc{0}{0}=0\in \mathbb{R}^{n_x}= \mathcal{E}_0$, and thus $\ktrigx\geq 1$. 
%Define the prediction error
%\begin{align}
%\label{TM:error}
%	\ermpc{j} = \trajmpc{w}{j} - \trajmpc{0}{j},
%\end{align}
%indicating the mismatch between the perturbed system and the nominal one. 
The closed-loop dynamics is then, for all $t\in\set{Z}_{\geq 0}$,
%The closed-loop dynamics given the initial state~$\xi_0$ is then
%\begin{align}
%\label{eq:mpc_dynamics}
%x^+_i=Ax_i+B\Umpci{i}+w_i,~\forall i\in\Seq{0}{\ktrigx},
%\end{align}
%where $x_0=x$.
\begin{subequations}
\label{eq:closed_dyn}
\begin{align}
\label{eq:mpc_dynamics}
\xi_{t+1}& = A \xi_t +B \Umpcdyn{t-\tau_t}{\xi_{\tau_t}} + w_t,\\
%\end{align}
%\begin{align*}
%\tau_{t+1} = \left[(\tau_t+1) \Ind{\hypreC{\tau_t}}{\xi_t-\trajmpC{0}{\tau_t}(\xi_{t-\tau_t})}\right] \cdot \left[\Ind{\Seq{0}{N-2}}{\tau_t}\right].
%\end{align*}
%Notice that $\tau_t\in\Seq{0}{N-1}$ is simply the inter-event, time counter. Consider $\Ttrig$ represents the triggering instance,
%\begin{align}
\label{eq:trig_evo}
\tau_{t+1}&= 
\left\lbrace
\begin{array}{lc}
\tau_t, & 
t-\tau_t\leq N-1 ~\text{and}~ \xi_t - \trajmpC{0}{t-\tau_t}(\xi_{\tau_t})\in \hypreC{t-\tau_t}(\xi_{\tau_t}), \\
%&\tau_t, & \big(t-\tau_t<N-1\big) \wedge \big(\xi_t - \trajmpC{0}{t-\tau_t}(\xi_{\tau_t})\in \hypreC{t-\tau_t}(\xi_{\tau_t})\big),\\
%&\tau_t, & \big(t-\tau_t<N-1\big) \wedge \big(\xi_t - \trajmpC{0}{t-\tau_t}(\xi_{t-\tau_t})\in \hypreC{t-\tau_t}\big),\\
t, & \text{otherwise},
%&t, & \big(t-\tau_t\leq N-1\big) \wedge \big(\xi_t - \trajmpC{0}{t-\tau_t}(\xi_{t-\tau_t})\notin \hypreC{t-\tau_t}\big),\\
%&t+1, & \big(t-\tau_t = N-1\big) \wedge \big(\xi_t - \trajmpC{0}{t-\tau_t}(\xi_{t-\tau_t})\in \hypreC{t-\tau_t}\big),
\end{array}
\right.
\end{align} 
\end{subequations}
given the initial state~$\xi_0$ and the initial triggering instance~$\tau_0=0$. Here, $\tau_t$ denotes the last triggering instance up to time~$t$. Also, notice that a mandatory triggering is put in place at time $\tau_{t}+N$. 
The problem addressed in this paper is now introduced. 
%\vspace{-.2cm}
\begin{Prob}
\label{prb:etMPC}
Consider the closed-loop dynamics~\eqref{eq:closed_dyn} under Assumptions~\ref{as:ass_0}-\ref{as:terminal}. 
%Consider the dynamics~\eqref{eq:dynamics} with the constraints~\eqref{eq:constraints}.
%Under Assumptions~\ref{as:ass_0}-\ref{as:terminal}, 
% now the sensory units are decentralized (in the sense that there is no centralized unit outside the controller unit that collects all the states).
%the controller synthesizes the optimal sequence~$\UmpC(\xi_{\tau_t})$ by solving the 
%the controller synthesizes the optimal sequence~$\Umpc$ by solving the 
%problem~\eqref{MPC_tightened_target}. 
Devise an approach to construct the sequence of triggering sets~$\HypreC(\xi_{\tau_t})$ in \eqref{eq:trig_state} such that the trajectories of the closed-loop dynamics satisfy:
\begin{itemize}
\item \textbf{\emph{Recursive feasibility:}} If $V^*_N(\xi_0)<\infty$, then $V^*_N(\xi_t)<\infty$, for all $t\in\mathbb{Z}_{\geq0}$;
%For all $i\in\Seq{1}{\ktrig(\xi_{\tau_t})+1}$, $\text{dom}\Big(V^*_N\big(\trajmpc{w}{i}\big) \Big) \neq \varnothing$; 
\item \textbf{\emph{Robust stability:}} The states and inputs of the closed-loop dynamics converge to the target sets $\Targ{X}$ and $\Targ{U}$, respectively ($\lim_{t\to \infty}V^*_N(\xi_{t})=0$). %(that is, the states of \eqref{eq:mpc_dynamics} converge to the state target set~$\Targ{X}$ and the control inputs converge to the input target set~$\Targ{U}$).
\end{itemize}
\end{Prob}

\begin{Rem}[Smart actuators and sensors]
The actuator and sensor units are ``{smart}" in the following sense. 
The actuator (sensor) units can buffer the time-stamped and packetized sequence $\UmpC(\xi_{\tau_t})$ ($\{ \trajmpC{0}{s}(\xi_{\tau_t})\oplus \hypreC{s}(\xi_{\tau_t}) \}_{s=1}^{N-1}$). 
%Each actuator unit receives and stores the time-stamped and packetized sequence of its corresponding input action in the sequence~$\UmpC(\xi_{\tau_t})$ from the controller. Then, 
The actuator units consecutively apply the input action~$\UmpcI{s-\tau_t}(\xi_{\tau_t})$ on the plant at each time~$s\in\Seq{\tau_t}{\tau_{t+1}-1}$. 
The sensor units evaluate the triggering condition $$\xi_s \notin \trajmpC{0}{s-\tau_t}(\xi_{\tau_t})\oplus \hypreC{s-\tau_t}(\xi_{\tau_t}),$$ at each time~$s\in\Seq{\tau_t+1}{\tau_{t}+N-1}$. When the triggering condition holds at some time~$s$, the sensors send the most recent states~$\xi_s$ to the controller and the triggering instance is set to $\tau_{t+1}=s$.
\end{Rem}

\begin{Rem}[Iteration Complexity]
RMPC problems with linear dynamics, a quadratic cost function, and polytopic constraints are \emph{quadratic programs} for which
%(the cost function is quadratic in decision variables and constraints are given by linear equalities or inequalities in terms of decision variables). Interestingly, 
dedicated solvers provide the complexity per iteration $\mathcal{O}(N(n_x+n_u)^3)$ 
% by exploiting the structure of the corresponding QP 
\cite{wang2010fast}. 
\end{Rem}

\section{Main Results}
\label{sec:eb_RMPC}
In this section, we provide several approaches to construct the sequence of sets~$\Hyprec$ which meets the requirements of Problem~\ref{prb:etMPC}.
%In this section, we provide the event-based implementation policy of the RMPC method~\eqref{MPC_tightened_target}.
To this end, we begin with describing a certain type of constrained optimization problem that produces $\Hyprec$. Based on these constructed sets, we then state the main theoretical results of this paper. 

\subsection{Construction of Hyper-Rectangles }
\label{sebsec:hyp_rec}
Let $j\in\Seq{1}{N-1}$. 
The procedure to construct each hyper-rectangle $\hyprec{j}$ comprises the parametric representation of $\hyprec{j}$, the definition of auxiliary quantities associated with $\hyprec{j}$, and finally the optimization problem to find $\hyprec{j}$. 

Notice that one way to represent a hyper-rectangle~$\hyprec{j}$ is
\begin{align*}
%\label{eq:hyper_rec}
	\hyprec{j} := \big\{ \epsilon \in \mathbb{R}^{n_x}:  -\ul{e}_{j}(x) \leq \epsilon \leq \ol{e}_{j}(x)  \big\},
\end{align*}
for some vectors~$\ul{e}_{j}(x),\ol{e}_{j}(x)\in\mathbb{R}^{n_x}_{\geq0}$ . In other words, each hyper-rectangle $\hyprec{j}$ is parameterized by $2n_x$ entries of $\ul{e}_j(x)$ and $\ol{e}_j(x)$.

Let us now introduce the auxiliary quantities involved in the derivation of $\hyprec{j}$. 
% Recall that $\mathbf{U}^*_{k|k}=\{u^*_{k+i|k}\}_{i=0}^{N-1}$ and $\mathbf{X}^*_{k|k}=\{x^*_{k+i|k}\}_{i=0}^{N}$ are available from solving $\mathcal{P}(x_k)$.
Let $\Acl:=(A+BF)$ be the nominal, closed-loop state matrix. Define the input sequence~$\Utilde$ and   
the associated state sequence~$\Trajcan$ as
\begin{subequations}
\begin{align}
\label{eq:extend_statefb1}
\Utildei{i}:=\left\lbrace
\begin{aligned}
&\Umpci{j+i}, & i\in\Seq{0}{N-j-1},\\
&F \Acl^{j+i-N} \trajmpc{0}{N}, & i\in\Seq{N-j}{N-1},
\end{aligned}
\right.
\end{align}
\begin{align}
\label{eq:extend_statefb2}
\Trajcani{i}:=\left\lbrace
\begin{aligned}
&\trajmpc{0}{j+i}, & i\in\Seq{0}{N-j},\\
& \Acl^{j+i-N} \trajmpc{0}{N}, & i\in\Seq{N-j+1}{N}.
\end{aligned}
\right.
\end{align}
Notice that the above \emph{candidate input sequence} is constructed by concatenating the last $N-j$ elements of $\Umpc$ with the nominal feedback~$F$ (recursively) applied to the optimal terminal state~$\trajmpc{0}{N}$.
\end{subequations}

Define $\Targprod{U} :=\prod_{i=0}^{N-1}\targ{U}{i}$ and $\Targprod{X}:=\prod_{i=0}^{N-1}\targ{X}{i}$. 
Denote now the projections of optimal state and input sequences $\Trajmpc{0}$ and $\Umpc$ onto their corresponding target sets by $\Slak{X} \in \Targprod{X}$ and $\Slak{U} \in \Targprod{U}$, where for all $i\in\Seq{0}{N-1}$,
\begin{align*}
\Slaki{X}{i}:=\Proj{\targ{X}{i}}{\trajmpc{0}{i}}{Q}, ~ 
\Slaki{U}{i}:=\Proj{\targ{U}{i}}{\Umpci{i}}{R}.
\end{align*}
Based on the above definition, the next two auxiliary quantities are defined as follows. Let $\Slaktilde{U}$ and $\Slaktilde{X}$ represent the projection of $\Utilde$ and $\Trajcan$ onto $\mathscr{T}^{\mathbb{U}}_N$ and $\mathscr{T}^{\mathbb{X}}_N$, respectively. We have 
\begin{subequations}
\label{eq:Slack}
\begin{align}
\label{eq:Slack_1}
\Slaktildei{U}{i} :=\left\lbrace
\begin{aligned}
& \Slaki{U}{j+i} , &\quad i\in\Seq{0}{N-j-1},\\
& \Utildei{i+j}, &\quad i\in\Seq{N-j}{N-1},
\end{aligned}
\right.
\end{align}
\begin{align}
\label{eq:Slack_2}
\Slaktildei{X}{i}	:=\left\lbrace
\begin{aligned}
& \Slaki{X}{j+i} , & i\in\Seq{0}{N-j},\\
& \Trajcani{j+i}, & i\in\Seq{N-j+1}{N-1}.
\end{aligned}
\right.
\end{align}
\end{subequations}
%\textbf{(1) Construction of optimal input and state trajectories:} We now aim to construct the following (with some abuse of notation) optimal trajectories
%$\mathbf{U}^*_{k+j|k}:=\{u^*_{k+j+i|k}\}_{i=0}^{N-1}$, and $\mathbf{X}^*_{k+j|k}:=\{x^*_{k+j+i|k}\}_{i=0}^{N}$.
%The term abuse of notation refers to the fact that we have only access to $u^*_{k+j+i|k}$ for $j+i\leq N-1$ and $x^*_{k+j+i|k}$ for $j+i\leq N$ from solving $\mathcal{P}(x_k)$. We adopt the convention
%\begin{subequations}
%\label{eq:extend_statefb}
%	\begin{align}
%	u^*_{k+j+i|k} &= F (A + BF)^{ j+ i - N} x^*_{k+N|k},&~\text{for}~ j+i\geq N,\\
%	x^*_{k+j+i|k} &= (A + BF)^{j + i - N} x^*_{k+N|k},&~\text{for}~ j+i\geq N+1,
%	\end{align}
%\end{subequations}
%to construct unavailable ``{optimal}" input and state trajectories based on the solution of $\mathcal{P}(x_k)$.
Let us clarify the conventions used in \eqref{eq:Slack}. 
Notice that the definition of $\Utildei{i}$ in \eqref{eq:extend_statefb1} implies that $\Utildei{i}\in \targ{U}{N-1} \subseteq \targ{U}{i}$, for all $i\in\Seq{N-j}{N-1}$. 
That is, the distance~$d_{R}(\Utildei{i}, \targ{U}{i})=0$, and hence, $\Proj{\targ{U}{i}}{\Utildei{i}}{R}=\Utildei{i}$, as given in \eqref{eq:Slack_1}. A similar line of reasoning has been used in \eqref{eq:Slack_2}. 
%The convention introduced in \eqref{eq:extend_statefb} implies that $u^*_{k+j+i|k}\in\mathcal{T}_{u,i}$ and $x^*_{k+j+i|k}\in\mathcal{T}_{x,i}$, for $j+i > N-1$ (since they are state-feedback extensions of the terminal state $x^*_{k+N|k}$). In light of this fact, $d(x^*_{k + j+i | k}, \mathcal{T}_{x, i}, Q)=d(u^*_{k + j+i | k}, \mathcal{T}_{u, i}, R)=0$ and hence, the choice made in the definition~\eqref{eq:extend_slack_statefb} becomes apparent.

We next adopt the feedback gains $\tilde{K}_i$ and the state-transition matrices $\tilde{L}_i$ defined as
%\footnote{Notice that $\tilde{L}_1 = A$: the error evolves in open-loop for one time step.}
\begin{subequations}
\label{eq:mat_tilde}
\begin{alignat}{3}
	&\tilde{K}_0 = 0_{n_u \times n_x}, ~\tilde{K}_{i+1} = K_{i},~ \forall i \in \Seq{0}{N-2} ,\label{K_tilde}\\
	&\tilde{L}_{0}=\mathsf{I}_{n_x},~\tilde{L}_{i+1} = (A + B\tilde{K}_i)\tilde{L}_{i},~ \forall i \in \Seq{0}{N-1}.\label{L_tilde}
\end{alignat}
\end{subequations}
In the following, we use the matrices~\eqref{eq:mat_tilde} to identify certain sets around the optimal state sequence~$\Trajmpc{0}$. These sets in turn will be used to formulate recursive feasibility and robust stability for the event-triggering setting (see the problem~\eqref{prob:decentralized_TM} and Section~\ref{subsec:th_rs}).

Let us now provide two definitions for the volume of $\hyprec{j}$, that are
\begin{subequations}
\label{eq:vol_def}
\begin{align}
\text{vol}_1(\hyprec{j})&:=\underset{p\in\Seq{1}{n_x}}{\prod}~\big(\ol{e}^p_{j}(x)+\ul{e}^p_{j}(x)\big), \label{eq:vol_def_1}\\
\text{vol}_2(\hyprec{j})&:=\underset{p\in\Seq{1}{n_x}}{\prod}~\big(\ol{e}^p_{j}(x) \times \ul{e}^p_{j}(x)\big), \label{eq:vol_def_2}
\end{align}
where $\ol{e}^p_{j}(x)$ (resp. $\ul{e}^p_{j}(x)$) denotes the $p$-th entry of $\ol{e}_{j}(x)$ (resp. $\ul{e}_{j}(x)$).
\end{subequations}
Notice that \eqref{eq:vol_def_1} is the standard definition of volume for $\hyprec{j}$ in $\mathbb{R}^{n_x}$. 
As it will be discussed later on, the application of \eqref{eq:vol_def_1} to construct $\hyprec{j}$ leads to a more asymmetric spread of $ \hyprec{j}$ around $\trajmpc{0}{j}$ compared to the application of \eqref{eq:vol_def_2}. The asymmetry in turn implies that the triggering mechanism has no robustness in certain error directions, see Remark~\ref{rem:dir_sens} for further details. Nonetheless, the definition~\eqref{eq:vol_def_1} leads to the construction of sets that have the maximum possible volume, in particular, higher than the ones constructed based on \eqref{eq:vol_def_2}. 
%The non-standard definition~\eqref{eq:vol_def_2} of the volume is introduced to handle the asymmetry issue and to promote a more symmetric construction of $\mathcal{E}_{j,k}$ around the origin compared to the construction based on the definition~\eqref{eq:vol_def_1}. 

For all $j\in\Seq{1}{N-1}$, the problem to find each $\hyprec{j}$ is
\begin{subequations}
\label{prob:decentralized_TM}
	\begin{alignat}{2}
%	&\max_{\ol{e}^p_{j,k}, \ul{e}^p_{j,k}}~ \left(\sum_{p=1}^{n} \ln( \ol{e}^p_{j,k} +\ul{e}^p_{j,k}) \right)\\
	&\max_{\ol{e}_{j}(x), \ul{e}_{j}(x)\geq 0}~ \text{vol}_q(\hyprec{j})\\
	&\text{s.t.}& \nonumber \\
%	& \ul{e}_{j} \geq 0, \bar{e}_{j} \geq 0, & \label{inclusion_error} \\
	& \Trajcani{i} \in \mathcal{X}_i \ominus \tilde{L}_i \hyprec{j}, & ~\forall i \in \Seq{0}{N-1}\label{inclusion_state}, \\
	& \Utildei{i} \in \mathcal{U}_i \ominus \tilde{K}_i \tilde{L}_i \hyprec{j}, & ~\forall i \in \Seq{0}{N-1} \label{inclusion_input}, \\
	&\Slaktildei{X}{i} \in \targ{X}{i} \ominus \tilde{L}_i \hyprec{j} , &  ~\forall i \in \Seq{0}{N-1} \label{inclusion_slack_state}, \\
	&\Slaktildei{U}{i} \in \targ{U}{i} \ominus \tilde{K}_i\tilde{L}_i \hyprec{j} , & ~\forall i \in \Seq{0}{N-1}\label{inclusion_slack_input},
	\end{alignat}
\end{subequations}
where $q\in\{1,2\}$ determines which type of the volume definition in \eqref{eq:vol_def} is chosen. 
%Notice that all the required data to solve \eqref{prob:decentralized_TM} is available after solving \eqref{MPC_tightened_target}.
%Notice that in order to solve the problem~(\ref{prob:decentralized_TM}) all required data is already available at instant $k$. Hence, it can be solved for $j \in \mathsf{N}_{[N - 1]}$ directly after an {RMPC} triggering instant. 
%\begin{Rem}[Non-convexity and parametric-in-set constraints]
%\label{rem:non_conv_hyper}
Notice that the objective function $\text{vol}_q(\hyprec{j})$ is a nonlinear, non-convex function with a decision variable~$\hyprec{j}$. Hence, the problem~\eqref{prob:decentralized_TM} is difficult to solve. 
In the next subsection, we show that this problem remains practically solvable, in particular, the set-based constraints \eqref{inclusion_state}-\eqref{inclusion_slack_input} are effectively representable by linear inequalities (i.e., polytopic inequalities) such that (i)~the optimization problem~\eqref{prob:decentralized_TM} has a CP counterpart (in Theorem~\ref{theo:lp-nlp-Box_1}), and (ii)~the optimization problem~\eqref{prob:decentralized_TM} admits an LP relaxation (in Theorem~\ref{theo:lp-nlp-Box_2}).
%\end{Rem}

\subsection{Event-Based Implementation}
\label{subsec:main_res}
We first show that robust stability of the event-triggered, closed-loop dynamics~\eqref{eq:closed_dyn} is guaranteed, which in turn leads to recursive feasibility of the closed-loop system. The triggering mechanism~\eqref{eq:trig_state} is constructed by the approach proposed in \eqref{prob:decentralized_TM}. 
%In simple words, if the incurred prediction errors caused by the model mismatch between two consecutive triggering instants are inside the hyper-rectangles $\mathcal{E}_{k}$, then, the decentralized, event-based implementation of the RMPC is both robustly recursively feasible and robustly stable.
We next establish that the non-convex problem~\eqref{prob:decentralized_TM} to construct the hyper-rectangles $\Hyprec$ has a CP reformulation and an LP relaxation, and therefore can be efficiently solved in practice. 
%Consider the perturbed LTI dynamics~\eqref{eq:dynamics} subject to the constraints~\eqref{eq:constraints}. Suppose that the initial state $x\in\domval$.
%For all $i\in \Seq{0}{N-1}$, let us denote $ \trajmpc{w}{i} - \trajmpc{0}{i}$ by $\ermpc{i}$. 
%\begin{Thm}[Robust recursive feasibility]
%	\label{theo_3}
%For all $t\in \mathbb{Z}_{>0}$, the problem~\eqref{MPC_tightened_target} remains feasible, i.e., $V^*_N(\xi_t)<\infty$.
%For all $j\in \Seq{0}{\ktrigx}$, the state~$\trajmpc{w}{j}\in\mathbb{X}$. Moreover, there exists an admissible, input sequence ~$\bm{u} \in \Set{U}_N$ such that the associated, state sequence~$\Trajfeas{0}{\trajmpc{w}{\ktrigx}}\in\mathscr{X}_N$, i.e., $\text{dom}\Big(V^*_N\big(\trajmpc{w}{\ktrigx}\big) \Big)\neq \varnothing$.
%\end{Thm}

\begin{Thm}[Robust convergence]
	\label{theo_4}
Consider the closed-loop dynamics~\eqref{eq:closed_dyn}, and suppose that the initial state $\xi_0$ is feasible (i.e., $V^*_N(\xi_0)<\infty$). For all $s\in\Seq{\tau_t+1}{\tau_{t+1}}$, there exists an  input sequence~$\bm{u}\in\Set{U}_N$ such that
\begin{align} \label{eq:rs_theo}
& V^*_N (\xi_{\tau_{t+1}} )-V^*_N(\xi_{\tau_t})  \leq  
V_N\big(\xi_s, \bm{u} (\xi_{s}) \big)-V^*_N(\xi_{\tau_t}) 
 \leq - \Big(\sum_{k=0}^{s-\tau_t-1} \DisT{Q}{\trajmpC{0}{k}(\xi_{\tau_t})}{X}{k} +\DisT{R}{\UmpcI{k}(\xi_{\tau_t})}{U}{k}\Big).
\end{align}
In particular, the closed-loop dynamics~\eqref{eq:closed_dyn} is asymptotically stable, i.e., $\lim_{t\to \infty}V^*_N(\xi_{t})=0$. 
%for $j=\ktrigx$, 
%\begin{align*}
%V_N^*(\trajmpc{w}{j})-V_N^*(x)\leq0.
%\end{align*}	
%For all $j\in \Seq{1}{\ktrigx-1}$, there exists an  input sequence~$\bm{u}\in\Set{U}_N$ such that
%\begin{multline*}
%V_N\big(\traJ{0}{0}{\trajmpc{w}{j}}, \bm{u}\big)-V^*_N(x) \leq \\
%- \sum_{k=0}^{j-1} \DisT{Q}{\trajmpc{0}{k}}{X}{k} +\DisT{R}{\Umpci{k}}{U}{k},
%\end{multline*}
%where $\bm{u}_0=\Umpci{j}$. Moreover, for $j=\ktrigx$, 
%\begin{align*}
%V_N^*(\trajmpc{w}{j})-V_N^*(x)\leq0.
%\end{align*}
\end{Thm}

\begin{Rem}[Recursive feasibility]
Notice that the second inequality in \eqref{eq:rs_theo} implies that $V_N\big(\xi_s, \bm{u} (\xi_{s}) \big)<\infty$, for all $s\in\Seq{\tau_t+1}{\tau_{t+1}}$. In other words, the optimization problem~\eqref{MPC_tightened_target} remains feasible for all time~$t\in \mathbb{Z}_{>0}$.  
\end{Rem}

%{\color{red}\begin{Rem}[Difference with standard analyses] 
%In standard RMPC, it is only required to show that the following problem~$V^*_N\big(\trajmpC{w}{1}(\xi_{t})\big)$ remains feasible. 
%, i.e., $\text{dom}\big(V^*_N\big(\trajmpc{w}{1}\big) \big)\neq \varnothing$. 
%In event-based case, one should further ensure that the state~$\trajmpC{w}{t-\tau_{t}}(\xi_{\tau_t})$ remain feasible for all~$t\in\Seq{\tau_t+1}{\tau_{t+1}}$. 
%$\trajmpc{w}{j}\in\mathbb{X}$ for all $j\in \Seq{0}{\ktrigx}$, i.e.,
%\end{Rem}}
\begin{Rem}[Transmission protocol]
We assume that all sensor and actuator units are clock-synchronized. 
%(by doing so, one can effectively reduce the listening time and conserve energy). 
When the problem~\eqref{MPC_tightened_target} is solved, the controller node sends: (i) $\UmpC(\xi_{\tau_t})$ to the actuator nodes and 
(ii) each entry of $\trajmpC{0}{j}(\xi_{\tau_t}) - \ul{e}_{j}(\xi_{\tau_t})$ and $\trajmpC{0}{j}(\xi_{\tau_t}) + \ol{e}_{j}(\xi_{\tau_t})$ to the corresponding sensory nodes, for all $j\in\Seq{1}{N-1}$. Moreover, the $n_x$ sensor units declare a triggering instance to each other, through a cost-efficient short-range transmission. Then, all sensors declare their time-stamped, observed states to the controller. 
%We assume that all communication units on sensors, controller and actuators are clock-synchronized (by doing so, one can effectively reduce the listening time and conserve energy). At each instant $k$ that the RMPC problem~\eqref{MPC_tightened_target} is solved, the controller node sends: (i) $\mathbf{U}^*_{k|k}$ to the actuator nodes and (2) the state individual bounds $x^{*,p}_{k+j|k}-\ul{e}^p_{j,k}$ and $x^{*,p}_{k+j|k}+\ol{e}^p_{j,k}$ to the corresponding sensory nodes, for all $j\in\mathsf{N}_{[N-1]}$. We further assume that the $n_x$ sensor units declare a triggering instant to each other, through a cost-efficient short-range transmission. Each individual sensor then declares its observed state to the controller. 
\end{Rem}

The successful usage of the above results is conditioned upon the premise that there exist computationally tractable methods to construct the sets $\Hyprec$.
We now revisit problem~\eqref{prob:decentralized_TM} to show that such a premise is valid by providing two frameworks: one in a CP form and another one in an LP form.
In these frameworks, the parametric-in-set constraints~\eqref{inclusion_state}-\eqref{inclusion_slack_input} can be reformulated into a new set of linear inequalities in terms of the vertices of each set $\hyprec{j}$.
We shall call the polytope represented by the derived linear inequalities, the \emph{principal} polytope $\bar{\mathcal{S}}$.
Both frameworks try to find a maximum-volume hyper-rectangle $\hyprec{j}$ inscribed (or contained) in the principal polytope such that $0\in\hyprec{j}$. 
%This problem is closely related to a well-studied problem in the literature known as ``{inradius}" of a polytopic set with respect to the polytopal norm induced by a hyper-rectangle with fixed (related to the LP form) or variable (related to the CP form) edge ratios. ( See e.g., \cite{gritzmann1993computational} and \cite{bemporad2004inner} for a detailed discussion on such problems.) 
In the LP framework, we partly employ some results from \cite{bemporad2004inner}, see Section~\ref{subsec:box_proof} and avoid reiterating the proofs of  borrowed material. 
% (citing the exact results of \cite{bemporad2004inner} that have been used in our proof).
%{\color{red}$\mathbf{x \in \mathbb{R}^p, \quad\quad \mathcal{S}\subset \mathbb{R}^p, \quad\quad A_{\mathcal{S}}\in \mathbb{R}^{m\times p},\quad\quad b_{\mathcal{S}}\in \mathbb{R}^m, \quad\quad M\in\mathbb{R}^{p\times k}, \quad\quad \mathcal{B}(l,u)\subset \mathbb{R}^k  }$}
For notational convenience, let $ \xi\in \mathcal{S}\ominus M \mathcal{B}(l,u)$ represent a concatenated version of the constraint~\eqref{inclusion_state}-\eqref{inclusion_slack_input} where, in particular, $\mathcal{B}(l,u):=\hyprec{j}$. Hereafter, when we take volume (of a hyper-rectangle) as defined in \eqref{eq:vol_def} with index $q=1$ and $q=2$ referring to \eqref{eq:vol_def_1} and \eqref{eq:vol_def_2}, respectively. 

\begin{Thm}[Volume maximization - CP reformulation]
\label{theo:lp-nlp-Box_1}
Consider a vector $\xi\in\mathbb{R}^p$, a matrix $M\in\mathbb{R}^{p \times k}$, and a polytope $\mathcal{S}=\{s\in\mathbb{R}^p:~ A_{\mathcal{S}}s\leq b_{\mathcal{S}} \}$ containing the origin where $A_{\mathcal{S}}\in\mathbb{R}^{m\times p}$ and $b_{\mathcal{S}}\in\mathbb{R}^m$. 
The maximum volume hyper-rectangle $\mathcal{B}(l,u)\subset \mathbb{R}^k$ that contains the origin and satisfies $\xi \in \mathcal{S}\ominus M \mathcal{B}(l,u)$ is $\mathcal{B}(-\ul{v}^*,\ol{v}^*)$ where $\ul{v}^*$ and $\ol{v}^*$ are the optimal solutions of the problem
\begin{equation}
\label{eq:NLP_form}
\begin{aligned}
  \underset{\ul{v},\ol{v}}{\min} \quad &f_q\big(\ol{v},\ul{v}\big) \\
 \text{\emph{s.t.}}               \quad   & \langle w^i, [\ol{v}^\top ~ \ul{v}^\top]^\top  \rangle \leq    b_{i,\mathcal{S}} - a_{i,\mathcal{S}} \xi, \forall i\in\Seq{1}{m}, \\
 & \ol{v}\geq0,~\ul{v}\geq0,                     
\end{aligned}
\end{equation}
where for $q\in\{1,2\}$
\begin{subequations}
\label{eq:cp_cost}
\begin{align}
f_1\big(\ol{v},\ul{v}\big)&:= -\underset{j\in\Seq{1}{k}}{\Sigma}~\log\big(\ol{v}_j+\ul{v}_j),\label{eq:cp_cost_1}\\
f_2\big(\ol{v},\ul{v}\big)&:= -\underset{j\in\Seq{1}{k}}{\Sigma}~\log\big(\ol{v}_j\big)+\log\big(\ul{v}_j), \label{eq:cp_cost_2}
\end{align}
\end{subequations}
and for all $j\in\Seq{1}{k}$
\begin{subequations}
\label{eq:NLP_form_w}
\begin{align}
w^i_j&=\left\lbrace
\begin{aligned}
&\big(M^\top a^{\top}_{i,\mathcal{S}}\big)_j, &\text{\emph{if}~} \hat{w}^i_j=1,\\
&0, &\text{\emph{otherwise}},
\end{aligned}
\right.\\
w^i_{k+j}&=\left\lbrace
\begin{aligned}
&-\big(M^\top a^{\top}_{i,\mathcal{S}}\big)_j, &\text{\emph{if}~} \hat{w}^i_j=-1,\\
&0, &\text{\emph{otherwise}},
\end{aligned}
\right.
\end{align}
\end{subequations}
with $\hat{w}^i := \text{\emph{sign}}(M^\top a^{\top}_{i,\mathcal{S}})$,  for all $i\in\Seq{1}{m}$.
\end{Thm}

\begin{Thm}[Volume maximization - LP relaxation]
\label{theo:lp-nlp-Box_2}
Suppose the hypotheses in Theorem~\ref{theo:lp-nlp-Box_1} hold. 
\begin{itemize}
\item $\mathbf{(q=1)}$ The maximum volume $r$-constrained hyper-rectangle $\mathcal{B}(l,u)\subset \mathbb{R}^k$ that contains the origin and satisfies $ \xi \in \mathcal{S}\ominus M \mathcal{B}(l,u)$ is $\mathcal{B}(z^*,z^*+\lambda^*r )$ for which $z^*\in\mathbb{R}^k$ and $\lambda^*\in\mathbb{R}$ are the optimal solution of the problem
\begin{subequations}
\label{eq:LP1_all}
\begin{equation}
\label{eq:LP_form}
\begin{aligned}
  \underset{z,\lambda}{\max} \quad &\lambda\\
 \text{\emph{s.t.}}               \quad   & A_{\mathcal{S}}M z + (A_{\mathcal{S}}M)^+ r \lambda \leq b_{\mathcal{S}} - A_{\mathcal{S}} \xi   \\
 & z+\lambda r \geq 0, ~z\leq 0,
\end{aligned}
\end{equation}
where the $j$-th entry of $r$, $j\in \Seq{1}{k}$, is defined as
\begin{equation}
\label{eq:LP_rj}
\begin{aligned}
r_j(\bar{\mathcal{S}})& := & \underset{z,\omega}{\max} \quad & \omega\\
                 & &\text{\emph{s.t.}} \quad & A_{\mathcal{S}}M z \leq b_{\mathcal{S}} - A_{\mathcal{S}} \xi\\
                 & & \quad & A_{\mathcal{S}}M (z+\omega e_j) \leq b_{\mathcal{S}} - A_{\mathcal{S}} \xi\\
                 & &\quad  & z + \omega e_j \geq0,  ~z \leq 0,\\
\end{aligned}
\end{equation}
\end{subequations}
where $e_j\in\mathbb{R}^k$ is the unit vector in the $j$-th direction and the polytope $\bar{\mathcal{S}}$ is
\begin{align*}
\bar{\mathcal{S}}:=\{z\in\mathbb{R}^k:~ A_{\mathcal{S}} M z \leq b_{\mathcal{S}} - A_{\mathcal{S}} \xi  \}.
\end{align*} 

\item $\mathbf{(q=2)}$ The maximum volume $r$-constrained hyper-rectangle $\mathcal{B}(l,u)\subset \mathbb{R}^k$ that contains the origin and satisfies $ \xi \in \mathcal{S}\ominus M \mathcal{B}(l,u)$ is $\mathcal{B}(-\lambda^*r_1,\lambda^*r_2 )$ for which $\lambda^*\in\mathbb{R}$ is the optimal solution of the problem
\begin{subequations}
\label{eq:LP2_all}
\begin{equation}
\label{eq:LP_form_1}
\begin{aligned}
  \underset{\lambda}{\max} \quad &\lambda\\
 \text{\emph{s.t.}}               \quad   & (W)^+ r \lambda \leq B   ,
\end{aligned}
\end{equation}
where $r=\big(r_2^\top, r_1^\top\big)^\top$ and the $j$-th entry of $r$, $j\in \Seq{1}{2k}$, is defined as
\begin{equation}
\label{eq:LP_rj_1}
\begin{aligned}
r_j & := & \underset{\omega}{\max} \quad & \omega\\
                 & &\text{\emph{s.t.}} \quad &  W' (\omega e_j) \leq B',
\end{aligned}
\end{equation}
\end{subequations}
where $e_j\in\mathbb{R}^{2k}$ is the unit vector in the $j$-th direction, 
\begin{align*}
&W = \big(w^1, \cdots, w^m \big)^\top, 
&W' =  
\left(\begin{aligned}
 &\qquad W \\
& \begin{aligned}
&-\mathsf{I}_k &  0_{k\times 1}\\
&0_{k\times 1}\;   &-\mathsf{I}_k\; 
\end{aligned} 
\end{aligned}\right), \\
&B = b_{\mathcal{S}} - A_{\mathcal{S}} \xi, 
&B'  = \big(B^\top, 0_{1\times 2k} \big)^\top,
\end{align*}
and for all $i\in\Seq{1}{m}$, $w^i$ are defined in \eqref{eq:NLP_form_w}.
\end{itemize}
\end{Thm}

We should emphasize that although Theorems~\ref{theo:lp-nlp-Box_1} \& \ref{theo:lp-nlp-Box_2} provide a way to construct $\hyprec{j}$ with a maximal volume, the derived set is not unique (the corresponding cost functions of these approaches are not strictly convex to guarantee the uniqueness of solution). 
In the remainder of the paper, we denote the construction approach based on the CP~\eqref{eq:NLP_form} with $q=1$ and $q=2$ by $\text{CP}_1$ and $\text{CP}_2$, respectively. 
Furthermore, $\text{LP}_1$ represents the LP relaxation~\eqref{eq:LP1_all} of $\text{CP}_1$ and  $\text{LP}_2$ denotes the LP relaxation~\eqref{eq:LP2_all} of  $\text{CP}_2$.

\subsection{Further Comments on Complexity and Sensitivity}
In the rest of this section, we allude briefly to two important practical aspects of the proposed construction approaches and possible directions to improve them. 
First, since these approaches are implemented online, they require an extra computation step besides the computation of the optimal input sequence. 
Notice that fixed-thresholding approaches in the literature, for example \cite{lehmann2013event}, avoid this extra step by considering pre-defined triggering sets. 
We provide the arithmetic complexity of the proposed approaches to quantify the extra computational burden. 
To this end, we adopt the following notion of an oracle to represent the optimization problems in this paper.
Let $A \in\mathbb{R}^{n_c\times n_d}$, $b \in \mathbb{R}^{n_c}$, $c \in\mathbb{R}^{n_d}$, and $f:\mathbb{R}^{n_d}\rightarrow\mathbb{R}$ be a concave function. 
Also, let $\textbf{lp}(n_c, n_d)$ denote the oracle complexity for solving $\max_{\eta}\{c^\top \eta: A\eta \leq b\}$, and $\textbf{cp}(n_c, n_d)$ denote the oracle complexity for solving $\max_{\eta}\{ f(\eta): A\eta \leq b\}$. 

\begin{Rem}[Computational complexity] \label{reb:oracle_comp}
The oracle complexity of the \text{CP} reformulations~\eqref{eq:NLP_form} in Theorem~\ref{theo:lp-nlp-Box_1} is $\mathbf{cp}(m+2k,2k)$ and of the \text{LP} reformulations~\eqref{eq:LP1_all} and \eqref{eq:LP2_all} in Theorem~\ref{theo:lp-nlp-Box_2} are $\mathbf{lp}(m+2k,k+1)+k\times\mathbf{lp}(2m+2k,k+1)$ and $\mathbf{lp}(m,1)+2k\times\mathbf{lp}(m+2k,1)$, respectively. 
A possible remedy to circumvent these computations is to introduce a \emph{state-independent} triggering law, as opposed to the current state-dependent law~\eqref{eq:trig_state}. This extension would allow to compute the desired sets offline and only once.
\end{Rem}

\begin{figure}[t!]
	\centering
	\subfigure[Well-shaped polytope example, $\frac{r_{\text{c}}}{r_{\circ}}=1.0861$.]{\label{fig:well_poly}\includegraphics[width=70mm]{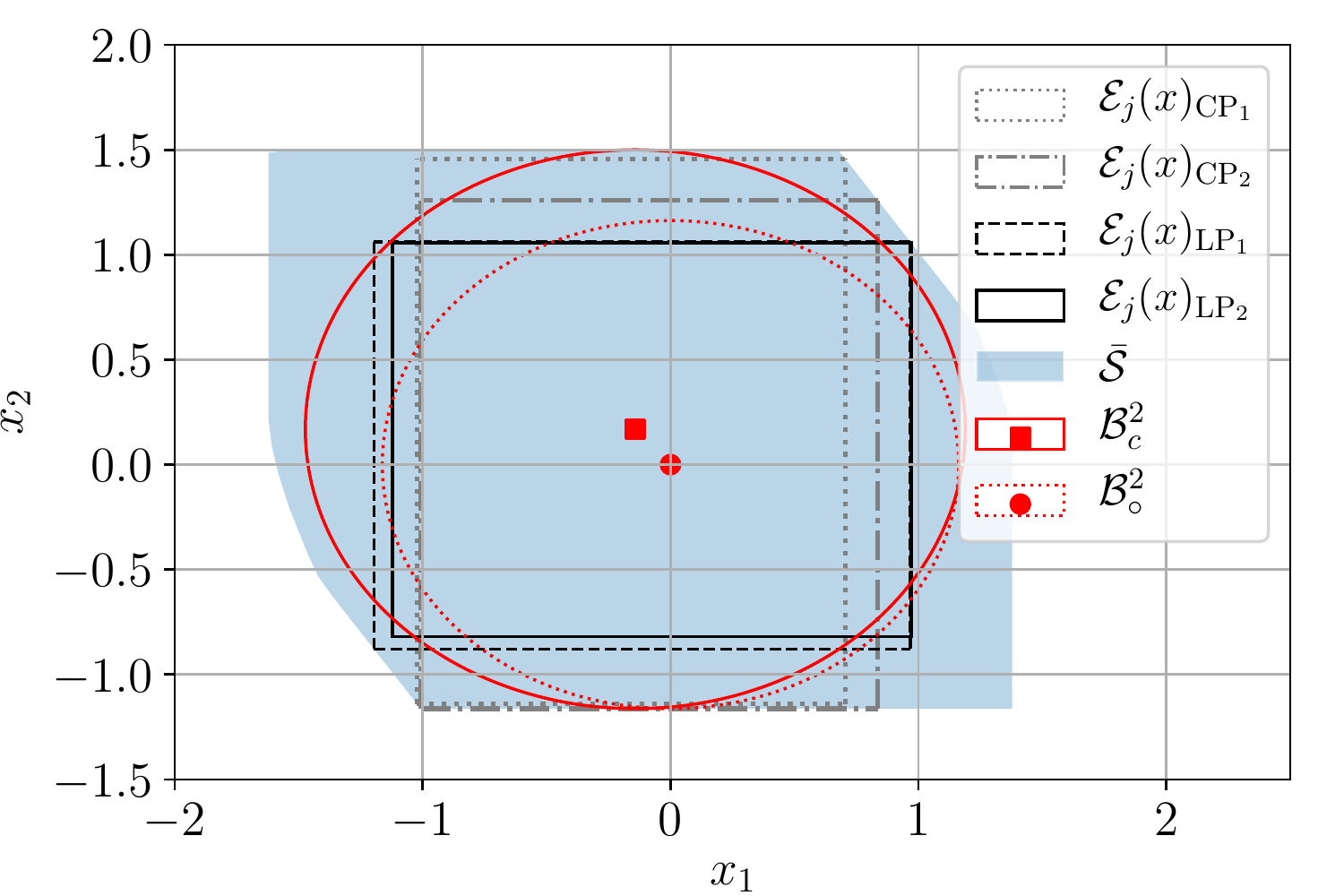}}    
%	\qquad
	\subfigure[Ill-shaped polytope example, $\frac{r_{\text{c}}}{r_{\circ}}=8.0669$.]{\label{fig:ill_poly}\includegraphics[width=70mm]{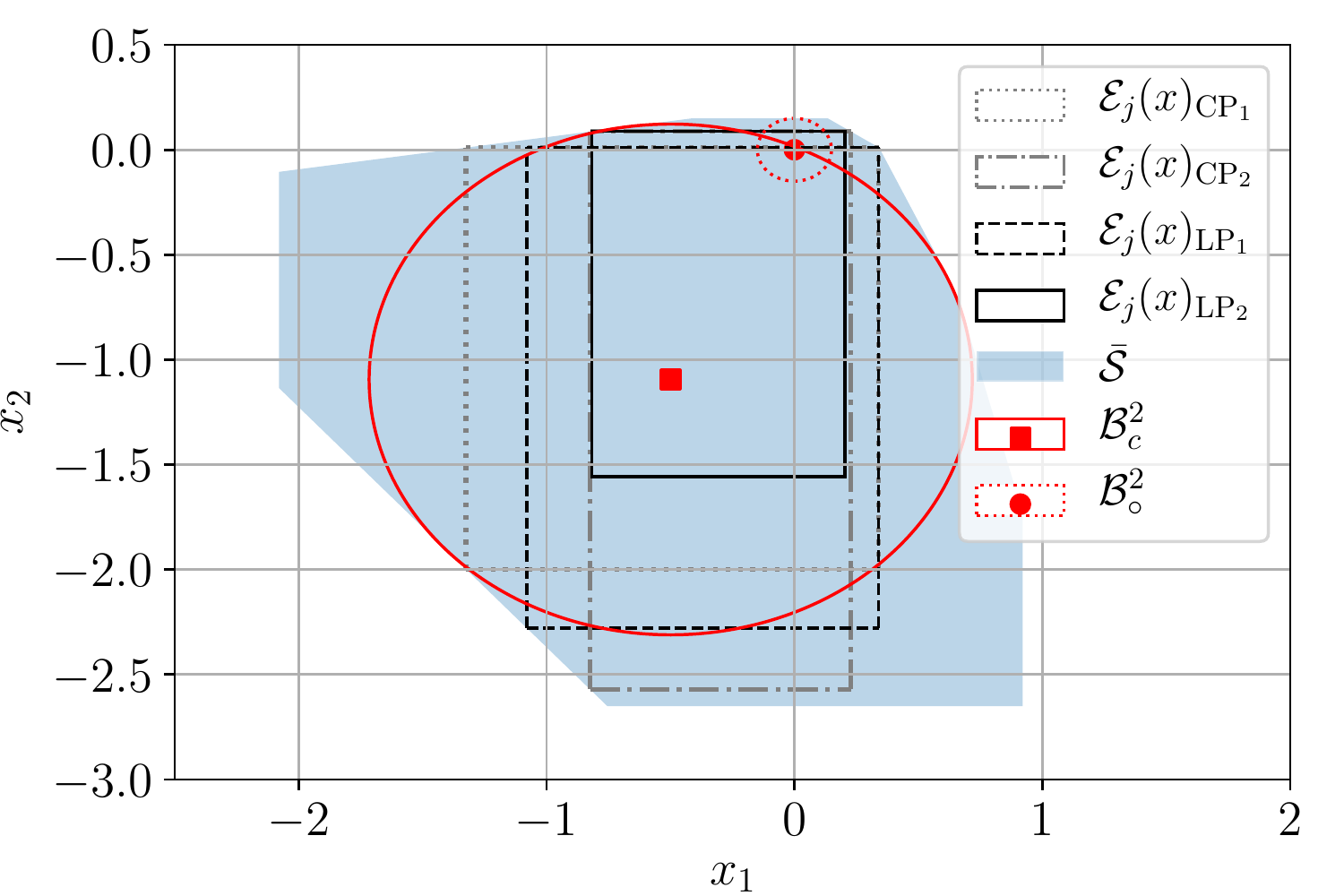}}
	\caption{Comparison of the CP and LP approaches to construct $\hyprec{j}\subseteq\ol{\mathcal{S}}$. \textbf{(a)} $\ol{\mathcal{S}}$ is distributed in a fairly uniform manner around the origin. All the approaches provide close behaviors. 
\textbf{(b)} $\ol{\mathcal{S}}$ is distributed in a relatively uneven manner around the origin. The approaches $\text{CP}_2$ and $\text{LP}_2$ promote more symmetric constructions compared to the approaches $\text{CP}_1$ and $\text{LP}_1$.}
	\label{fig:sensitivity}
\end{figure}

The other issue regarding the proposed approaches is the asymmetry of the triggering sets with respect to the optimal state sequence. 
Let polytope~$\ol{\mathcal{S}}\subset \mathbb{R}^{n_x}$ represent the constraints~\eqref{inclusion_state}-\eqref{inclusion_slack_input} that the triggering set~$\hyprec{j}$ satisfies. 
In other words, $\hyprec{j}$ is constructed inside $\ol{\mathcal{S}}$. 
Recall that $\hyprec{j}$ represents the ``{allowable}" prediction error so that the triggering mechanism is not activated. 
Qualitatively speaking, for a ``{better}" directional resilience against prediction errors, one would prefer symmetry in the constructed~$\hyprec{j}$. The above statements are schematically depicted in Figure~\ref{fig:sensitivity}. 
When $\ol{\mathcal{S}}$ is well-shaped as in Figure~\ref{fig:well_poly}, the approaches in Theorems~\ref{theo:lp-nlp-Box_1} \& \ref{theo:lp-nlp-Box_2} lead to a relatively symmetric set~$\hyprec{j}$ with respect to the origin. 
%As a result, the triggering mechanism shows the same level of sensitivity to opposite directions of the prediction error along different coordinates. 
When $\ol{\mathcal{S}}$ is ill-shaped as in Figure~\ref{fig:ill_poly}, the constructed set~$\hyprec{j}$ is however extremely asymmetric with respect to the origin along some coordinates. 
%This in turn leads to higher sensitivity of the triggering mechanism to a particular direction of the prediction error in those coordinates. 
This difference is well-captured by the geometric measure~$\frac{r_{\text{c}}}{r_{\circ}}$ of $\ol{\mathcal{S}}$, where $r_{\text{c}}$ is the radius of the maximal $2$-norm ball inside $\ol{\mathcal{S}}$, and $r_{\circ}$ is the radius of the maximal $2$-norm ball, centered at the origin and inside $\ol{\mathcal{S}}$. 
By definition, we have $\frac{r_{\text{c}}}{r_{\circ}}\geq 1$. 
Observe that in well-shaped cases $r_c/r_{\circ} \approx 1$ and in ill-shaped cases $r_c/r_{\circ} \gg 1$.

\begin{Rem}[Directional sensitivity to prediction errors]
\label{rem:dir_sens}
The directional sensitivity issue is the main reason for introducing the second definition~\eqref{eq:vol_def_2} of the volume. 
To see this, assume first that our goal is to maximize the $\log$ value of the volume of $\hyprec{j}$.  
Notice that the first definition~\eqref{eq:vol_def_1} solely aims for maximizing the width of $\hyprec{j}$ within $\ol{\mathcal{S}}$ along each coordinate. 
On the other hand, the second definition~\eqref{eq:vol_def_2} maximizes the width of $\hyprec{j}$ in both positive and negative directions along each coordinate.   
As shown in Figure~\ref{fig:sensitivity}, in both cases the set~$\hyprec{j}$ constructed by the approaches~$\text{CP}_2$ and $\text{LP}_2$ is typically more symmetric compared to those constructed by the approaches $\text{CP}_1$ and $\text{LP}_1$. 
An interesting research direction to alleviate this sensitivity issue is to investigate the impact of the MPC design parameters (e.g., the tightening gains~$\bm{K}$ or the target sets~$\Targ{X}$ and $\Targ{U}$).
\end{Rem}

\section{Technical Proofs}
\label{sec:proof}

\subsection{Proof of Theorem~\ref{theo_4}}
\label{subsec:th_rs}
%\begin{proof}

The proof consists of five main steps. Each step is labeled by the guaranteed property. 
Let $x:=\xi_{\tau_t}$ be the state at the last triggering instance.
Define the prediction error
\begin{align}
\label{TM:error}
	\ermpc{j} = \trajmpc{w}{j} - \trajmpc{0}{j},
\end{align}
indicating the mismatch between the perturbed system and the nominal one. 
For some integer~$j\in\Seq{0}{N-1}$, suppose that the mechanism is enabled at time~$j+1$, that is, either (1)~$j<N-1$ so that for all $i\in \Seq{0}{j}$, $\ermpc{i}\in \hyprec{i}$ and $\ermpc{j+1}\notin \hyprec{j+1}$, or (2)~$j=N-1$ (see equation~\eqref{eq:trig_evo}). 
%In what follows, we only consider the first case since a similar line of arguments for the second case can be used. 
%Recall that for $j=N-1$, the triggering mechanism is activated regardless of the state of the system. 
We omit the arguments of variables for convenience when it is clear from the context (unless mentioned otherwise). In what follows, we also use the notation $\ell(x_i,u_i)$ for $\DisT{Q}{x_i}{X}{i}+\DisT{R}{u_i}{U}{i} $ for notational simplicity.  

%Let $j\in\Seq{1}{N-2}$. Suppose that the mechanism is enabled at time~$j+1$, i.e., for all $i\in \Seq{1}{j}$, $\ermpc{i}\in \hyprec{i}$ and $\ermpc{j+1}\notin \hyprec{j+1}$. 
%Moreover, let the disturbance sequence~$\bm{w'}$ be $\{\ermpC{j}, 0, \ldots, 0\}$.

\textbf{1)~Inter-event recursive feasibility:} Define the \emph{candidate} input sequence~$\Ufeas$ such that for all $i\in\Seq{0}{N-1}$,
\begin{subequations}
\label{eq:xu_can}
\begin{align}
\label{eq:u_can}
\UfeasI{i} :=\UtildeI{i} + \tilde{K}_i\tilde{L}_i \ermpC{j},
\end{align}
and its associated \emph{candidate} state sequence~$\Xfeas$, where for all $i\in\Seq{0}{N}$,
\begin{align}
\label{eq:x_can}
\XfeasI{i}:= \TrajcanI{i} + \tilde{L}_i \ermpC{j}.
\end{align}
\end{subequations}
Note that $\XfeasI{0}=\trajmpC{w}{j}$ and $\UfeasI{0}=\UmpcI{j}$. 
We now establish that the sequences $\UfeaS$ and $\XfeaS$ satisfy $\UfeaS\in\Set{U}_N$ and $\XfeaS\in\Set{X}_N$, i.e., $V_N(\trajmpC{w}{j},\UfeaS)<\infty$. 
%{\color{blue}($\Xfeasi{0}$ and $\Ufeasi{0}$ are the actual observed state and applied input of the closed-loop dynamics, respectively).} 
By assumption, $\ermpC{j} \in \hypreC{j}$. Moreover, $\hypreC{j}$ satisfies \eqref{inclusion_state}-\eqref{inclusion_input}. 
From the definition of the Pontryagin difference, it follows that $\UfeasI{i} \in \mathcal{U}_i$ and $\XfeasI{i}\in \mathcal{X}_i$, for all $i\in \Seq{0}{N-1}$. 
%From Assumption~\ref{as:terminal},  we have $\Xfeasi{N-1}\in \mathcal{X}_f$. 
Recall that $L_{N-1}=0$. Hence, $\tilde{L}_N=0$ and $\XfeasI{N}=\TrajcanI{N}$. From \eqref{eq:extend_statefb2}, we have $\TrajcanI{N}=\Acl^{j} \trajmpc{0}{N}$. 
Since $\mathcal{X}_f$ is a control invariant set, $\XfeasI{N}\in \mathcal{X}_f$. 
We conclude that $V_N(\trajmpC{w}{j},\UfeaS)<\infty$.

\textbf{2)~Inter-event cost function decay:} 
Observe that
\begin{align*}
\DisT{Q}{\XfeasI{i}}{X}{i} =\DisT{Q}{\TrajcanI{i} + \tilde{L}_i \ermpC{j}}{X}{i}  \leq \DisTM{Q}{\TrajcanI{i}}{X}{i},
\end{align*}
where we made use of the definition~\eqref{eq:x_can} and Lemma~\ref{lem:lem_1}, respectively. 
Recall from \eqref{inclusion_slack_state} that $\SlaktildeI{X}{i} \in \targ{X}{i} \ominus \tilde{L}_i \hypreC{j}$. 
Hence, 
\begin{subequations}
\label{eq:pr}
\begin{align}
\label{eq:pr1}
0 \leq \DisT{Q}{\XfeasI{i}}{X}{i} \leq d_{Q}(\TrajcanI{i} , \SlaktildeI{X}{i}).
\end{align}
Similarly, one can arrive at
\begin{align}
\label{eq:pr2}
0 \leq \DisT{R}{\UfeasI{i}}{U}{i} \leq d_{R}(\UtildeI{i} , \SlaktildeI{U}{i}).
\end{align}
\end{subequations}

Consider $i\in\Seq{0}{N-j-1}$. In light of the definitions~\eqref{eq:extend_statefb2} and \eqref{eq:Slack_2}, 
we have $d_{Q}(\TrajcanI{i} , \SlaktildeI{X}{i})=\DisT{Q}{\trajmpC{0}{j+i}}{X}{j+i}$. Thus,
\begin{subequations}
\label{eq:val_feas}
\begin{align}
\label{eq:val_feas1}
\DisT{Q}{\XfeasI{i}}{X}{i} \leq \DisT{Q}{\trajmpC{w}{j+i}}{X}{j+i}.
\end{align}
In a similar fashion, we can show
\begin{align}
\label{eq:val_feas2}
\DisT{R}{\UfeasI{i}}{U}{i} \leq \DisT{Q}{\UmpcI{j+i}}{U}{j+i}.
\end{align}
\end{subequations}
From \eqref{eq:val_feas}, it is then straightforward that
\begin{subequations}
\label{eq:val_can}
\begin{align}
\label{eq:val_can1}
\sum_{i=0}^{N-j-1} \ell(\XfeasI{i},\UfeasI{i}) \leq \sum_{i=0}^{N-j-1} \ell(\trajmpC{0}{j+i},\UmpcI{j+i}).
\end{align}

Now, let $i\in\Seq{N-j}{N-1}$ and consider the definition~\eqref{eq:Slack}. Then, $d_{Q}(\TrajcanI{i} , \SlaktildeI{X}{i})=d_{R}(\UtildeI{i} , \SlaktildeI{U}{i})=0$. These equality relations coupled with \eqref{eq:pr} give rise to 
\begin{align}
\label{eq:val_can2}
\sum_{i=N-j}^{N-1} \ell(\XfeasI{i},\UfeasI{i})=0.
\end{align}
\end{subequations}
From \eqref{eq:val_can}, we finally infer that if $\ermpC{j}\in \hypreC{j}$, then
\begin{align}
\label{eq:ineq_thm_1}
\Valf{\trajmpC{w}{j}}{\UfeaS} =\Valf{\XfeasI{0}}{\UfeaS}  \leq \Valopt{x} - \sum_{i=0}^{j-1} \ell(\trajmpC{0}{i},\UmpcI{i}).
\end{align}

\textbf{3)~At-event recursive feasibility:} Consider now the \emph{new candidate} input sequence~$\Uhcan$ where
\begin{subequations}
\label{eq:xu_canh}
\begin{align}
\label{eq:u_canh}
\Uhcani{i}:=
& \left\{
\begin{aligned}
&\UfeasI{i+1} +K_i L_i w_{j}, &  \forall i\in \Seq{0}{N-2},\\
&F \XfeasI{N}+FL_{N-1} w_j , &  i=N-1, %+ F  L_{N-1} w_{j}
\end{aligned}\right.
\end{align}
and its associated \emph{candidate} state sequence~$\Xhcan$ such that
\begin{align}
\label{eq:x_canh}
\XhcanI{i}:= \left\{
\begin{aligned}
&\XfeasI{i+1}+ L_i w_{j}, 
& \forall i\in \Seq{0}{N-1},\\
&\Acl \XhcanI{N-1} , & i=N,
\end{aligned}\right.
\end{align}
\end{subequations}
where $w_{j}\in\set{W}$. Observe that 
\begin{align*}
\XhcanI{0}&=\XfeasI{1}+w_{j}=A \XfeasI{0}+B\UfeasI{0}+w_{j}\\
&=A\trajmpC{w}{j}+B\UmpcI{j}+w_j=\trajmpC{w}{j+1}.
\end{align*}
We now show that $\UhcaN\in\Set{U}_N$ and $\XhcaN\in\Set{X}_N$, i.e., $V_N(\trajmpC{w}{j+1},\UhcaN)<\infty$. 
Observe that $\UfeasI{i+1}\in \mathcal{U}_{i+1}$ and $w_j\in\set{W}$. 
Hence, $\UfeasI{i+1}\in \mathcal{U}_{i+1}\oplus K_i L_i \set{W}$ for all $i\in\Seq{0}{N-2}$. 
Since $\mathcal{U}_{i+1}=\mathcal{U}_i\ominus K_i L_i \set{W}$, we have $\Uhcani{i}\in\mathcal{U}_i$ for all $i\in\Seq{0}{N-2}$. 
Recall now $\XfeasI{N}\in \mathcal{X}_f$ (from Step~1). Assumption~\ref{as:terminal} along with $L_{N-1}=0$ imply that $\Uhcani{N-1}\in \mathcal{U}_{N-1}$. 
We have $\XfeasI{i+1}\in \mathcal{X}_{i+1}$ for all $i\in\Seq{0}{N-2}$. Then, $\XhcanI{i}\in \mathcal{X}_{i+1} \oplus L_i \set{W}$. 
For all $i\in\Seq{0}{N-2}$, it follows from $\mathcal{X}_{i+1}=\mathcal{X}_i \ominus L_i \set{W}$ that $\XhcanI{i}\in\mathcal{X}_i$. 
Recall that $\XfeasI{N}\in\mathcal{X}_f$ and $L_{N-1}=0$. Hence, we arrive at $\XhcanI{N-1}\in \mathcal{X}_f$ and as a result $\XhcanI{N}\in\mathcal{X}_f$. 
We thus have $\UhcaN\in\Set{U}_N$ and $\XhcaN\in\Set{X}_N$, i.e.,  $V_N(\trajmpC{w}{j+1},\UhcaN)<\infty$.
%\end{proof}

%\subsection{Proof of Theorem~\ref{theo_4}}
%\label{subsec:th_rs}
%\begin{proof}
%Let $j\in\Seq{1}{N-2}$. Suppose that the mechanism is enabled at time~$j+1$. 

\textbf{4)~At-event value function decay:} Consider now $\UhcaN$ and $\XhcaN$ as the candidate input and state sequences at time~$j+1$, respectively. For all $i\in\Seq{0}{N-2}$ and for all $w_j\in\set{W}$,
\begin{subequations}
\label{eq:trig_val_dec}
\begin{align}
 \DisT{Q}{\XhcanI{i}}{X}{i} & = \DisT{Q}{\XfeasI{i+1}+L_i w_j}{X}{i} \label{eq:valtrig_1}\\
  & \leq \DisTM{Q}{\XfeasI{i+1}}{X}{i}= \DisT{Q}{\XfeasI{i+1}}{X}{i+1},\nonumber
\end{align}
where the first inequality follows from \eqref{eq:x_canh}, the inequality is implied by Lemma~\ref{lem:lem_1}, and the second equality is derived from \eqref{tight:state}. 
Following a similar argument, we arrive at
\begin{align}
\DisT{R}{\Uhcani{i}}{U}{i} \leq \DisT{R}{\UfeasI{i+1}}{U}{i+1}.
\end{align}
\end{subequations}
Since $L_{N-1}=0$, $\XhcanI{N-1}= \XfeasI{N}$ and $\Uhcani{N-1}=F \XfeasI{N}$. In Step~3, it is shown that $\XhcanI{N-1}\in\mathcal{X}_f$. Then Asuumption~\ref{as:terminal} implies that $\XhcanI{N-1}\in \targ{X}{N-1}$ and $\Uhcani{N-1}\in \targ{U}{N-1}$. Hence, $\ell(\XhcanI{N-1} ,\Uhcani{N-1}) =0$. By virtue of the inequalities in \eqref{eq:trig_val_dec}, we then arrive at 
\begin{align}
\label{eq:sum_val_dec}
\begin{aligned}
\Valf{\trajmpC{w}{j+1}}{\UhcaN} 
& =  \Valf{\XhcanI{0}}{\UhcaN}  \leq \sum_{i=1}^{N-1} \ell(\XfeasI{i},\UfeasI{i})\\
& = \Valf{\XfeasI{0}}{\UfeaS} - \ell(\XfeasI{0},\UfeasI{0})\\
& = \Valf{\XfeasI{0}}{\UfeaS} - \ell(\trajmpC{0}{j},\UmpcI{j}) \\
&  \leq \Valopt{x} - \sum_{i=0}^{j} \ell(\trajmpC{0}{i},\UmpcI{i}).
\end{aligned}
\end{align}
It follows from the optimality principle that $ \Valopt{\trajmpC{w}{j+1}} \leq \Valf{\trajmpC{w}{j+1}}{\UhcaN}$. This inequality along with \eqref{eq:sum_val_dec} in turn implies that 
\begin{align} 
\label{eq:ineq_thm_2}
\Valopt{\trajmpC{w}{j+1}} \leq\Valopt{x} - \sum_{i=0}^{j} \ell(\trajmpC{0}{i},\UmpcI{i}).
\end{align}

\textbf{5)~Robust convergence:} First, observe that \eqref{eq:rs_theo} is an immediate consequence of \eqref{eq:ineq_thm_1} and \eqref{eq:ineq_thm_2}. 
Let us now recall that $x=\xi_{\tau_t}$ and $\trajmpc{w}{j+1}=\xi_{\tau_{t+1}}$. 
Then, one can rewrite \eqref{eq:ineq_thm_2} as follows:
\begin{align*}
\Valopt{\xi_{\tau_{t+1}}}-\Valopt{\xi_{\tau_t}}\leq - \sum_{i=0}^{\tau_{t+1}-\tau_t-1} \ell\big(\trajmpC{0}{i}(\xi_{\tau}),\UmpcI{i}(\xi_{\tau})\big).
\end{align*}
Notice that the right-hand side of the above inequality is strictly negative unless when $\trajmpC{0}{i}(\xi_{\tau})\in\targ{X}{i}$ and $\UmpcI{i}(\xi_{\tau})\in \targ{U}{i}$ for all $i\in\Seq{0}{\tau_{t+1}-\tau_t-1}$. 
%Notice that the right side of the above inequality is, for all $j\in\Seq{0}{j}$: (1) negative when $\trajmpC{0}{i}(\xi_{\tau})\notin\targ{X}{i}$ and $\UmpcI{i}(\xi_{\tau})\notin \targ{U}{i}$, and (2) zero when $\trajmpC{0}{i}(\xi_{\tau})\in\targ{X}{i}$ and $\UmpcI{i}(\xi_{\tau})\in \targ{U}{i}$. 
Since $\Valopt{\xi_{\tau_{t+1}}}$ is a non-negative value, it is straightforward to observe that the states and inputs of the closed-loop dynamics~\eqref{eq:closed_dyn} converge to their corresponding target sets. This concludes the proof. 
%
%the amount of reduction in the cost function $J$ is $d(x_{k+j}, \mathcal{T}_{x, 0}, Q) + d(u^*_{k+j|k+j}, \mathcal{T}_{u, 0}, R)$. 
%Since the cost function $J\geq 0$, the above analysis guarantees that $J(x_{k_{\text{\emph{trig}}}}, \mathbf{U}^*_{k_{\text{\emph{trig}}}|k_{\text{\emph{trig}}}})$ converges to a steady value as the triggering instant $k_{\text{\emph{trig}}}\rightarrow \infty$. This in turn implies $d(x_{k+j}, \mathcal{T}_{x, 0}, Q) + d(u^*_{k+j|k+j}, \mathcal{T}_{u, 0}, R) \rightarrow 0$, as $k \rightarrow \infty$, from which we can infer that $x_k \rightarrow \mathbb{T}_x$ and $u_k \rightarrow \mathbb{T}_u$ as $k \rightarrow \infty$. This concludes the proof.
%\end{proof} 

\subsection{Proof of Theorems~\ref{theo:lp-nlp-Box_1} \& \ref{theo:lp-nlp-Box_2}}
\label{subsec:box_proof}
%\begin{proof}
We first begin with a preliminary argument that is shared between both theorems. We then carry on with the proof of each case in an orderly fashion. Notice that $\xi\in \mathcal{S}\ominus M \mathcal{B}(l,u)$ and $\mathcal{S}$ is a polytope by the theorems' hypothesis. 
%By virtue of the relation~\eqref{eq:pontDH}, 
Let $h_{M\mathcal{B}}$ be the support function of $M\mathcal{B}$. One can infer that   
\begin{align*}
\langle a_{i,\mathcal{S}}^\top , \xi \rangle & \leq b_{i,\mathcal{S}} - h_{M\mathcal{B}}(a_{i,\mathcal{S}}^\top), \forall i\in\Seq{1}{m}. 
\end{align*}
Next, observe that $\mathcal{B}(l,u)\subset \mathbb{R}^k$ is a polytope (and as a result bounded), and the domain $\mathcal{K}_{\mathcal{B}}$ on which the support function $h_{\mathcal{B}}$ is defined is the whole space, i.e., $\mathcal{K}_{\mathcal{B}}=\mathbb{R}^k$. Hence, $
h_{M\mathcal{B}}(a^\top_{i,\mathcal{S}}) = h_{\mathcal{B}}(M^\top a^\top_{i,\mathcal{S}})$,
and as a consequence
\begin{align*}
\langle a_{i,\mathcal{S}}^\top , \xi \rangle \leq b_{i,\mathcal{S}} - h_{\mathcal{B}}(M^\top a_{i,\mathcal{S}}^\top), \forall i\in\Seq{1}{m}.
\end{align*}
Rearranging the above inequality, we arrive at
\begin{align*}
h_{\mathcal{B}}(M^\top a_{i,\mathcal{S}}^\top)  \leq  b_{i,\mathcal{S}} - \langle a_{i,\mathcal{S}}^\top , \xi \rangle, \forall i\in\Seq{1}{m},
\end{align*}
where the only unknown entity is $h_{\mathcal{B}}(M^\top a_{i,\mathcal{S}}^\top)$ with $M^\top a_{i,\mathcal{S}}^\top \in \mathbb{R}^k$. It follows from the definition of the support function that $\langle M^\top a_{i,\mathcal{S}}^\top, z \rangle \leq h_{\mathcal{B}}(M^\top a_{i,\mathcal{S}}^\top) $ for all $z\in \mathbb{R}^k$. Thus,
\begin{align}
\label{eq:parBineq_1}
\langle M^\top a_{i,\mathcal{S}}^\top, z \rangle \leq    b_{i,\mathcal{S}} - \langle a_{i,\mathcal{S}}^\top , \xi \rangle, \forall i\in\Seq{1}{m},\forall z\in\mathcal{B}.
\end{align}
Let us now define for all $i\in\Seq{1}{m}$, $
a^{\top}_{i,\bar{\mathcal{S}}} := M^\top a^{\top}_{i,\mathcal{S}}$, $
b_{i,\bar{\mathcal{S}}}  := b_{i,\mathcal{S}} - \langle a^{\top}_{i,\mathcal{S}}, \xi  \rangle $, 
and the convex polytope (which we referred to as the \emph{principal} polytope in the paragraph before Theorem~\ref{theo:lp-nlp-Box_1})
\begin{equation}
\label{eq:parBineq_2}
\begin{aligned}
\bar{\mathcal{S}} &:= \{s\in \mathbb{R}^k:~ \langle a^{\top}_{i,\bar{\mathcal{S}}} , s  \rangle \leq b_{i,\bar{\mathcal{S}}}, \forall i\in \Seq{1}{m}  \}\\
 & \; = \{s\in \mathbb{R}^k:~ A_{\bar{\mathcal{S}}} s \leq b_{\bar{\mathcal{S}}}  \},
\end{aligned}
\end{equation}
where $A_{\bar{\mathcal{S}}}:=[a^{\top}_{1,\bar{\mathcal{S}}}, \cdots, a^{\top}_{m,\bar{\mathcal{S}}}]^\top=(M^\top A_{\mathcal{S}}^{\top})^{\top}=A_{\mathcal{S}}M$ and $b_{\bar{\mathcal{S}}}:=[b_{1,\bar{\mathcal{S}}}, \cdots,b_{m,\bar{\mathcal{S}}}]^\top= b_{\mathcal{S}} - A_{\mathcal{S}}\xi $. Now, one can deduce from the inequalities~\eqref{eq:parBineq_1} and the definition~\eqref{eq:parBineq_2} that the convex polytope $\bar{\mathcal{S}}$ contains the hyper-rectangle $\mathcal{B}(l,u)$, i.e., $\mathcal{B}(l,u)\subseteq \bar{\mathcal{S}}$. 
Notice that $\mathcal{B}(l,u)$ is parametric in the variables $l$ and $u$.

\textbf{Theorem~\ref{theo:lp-nlp-Box_1}:} In the CP framework, we propose a convex nonlinear program to compute the hyper-rectangle $\mathcal{B}(l,u)\subseteq \bar{\mathcal{S}}$ such that its volume is maximized. 
Suppose $\mathcal{B}(l,u)$ is parameterized as $l:=-\ul{v}=[-\ul{v}_1,\cdots,-\ul{v}_k]^\top$ and $u:=\ol{v}=[\ol{v}_1,\cdots,\ol{v}_k]^\top$ such that for all $i\in\Seq{1}{k}$, $\ul{v}_i$ and $\ol{v}_i$ are positive scalars (this condition has to do with the fact that the resulting hyper-rectangle should contain the origin). 
Recall the inequality~\eqref{eq:parBineq_1}, that is $\langle M^\top a_{i,\mathcal{S}}^\top, z \rangle \leq    b_{i,\mathcal{S}} - a_{i,\mathcal{S}} \xi$, for all $i\in\Seq{1}{m}$ and for all $z\in\mathcal{B}$. 
In what follows, we show that although the hyper-rectangle $\mathcal{B}(l,u)=\mathcal{B}(-\ul{v},\ol{v})$ is parametric, one can provide a closed form for its support function evaluated at $M^\top a^{\top}_{i,\mathcal{S}}$. By definition of a support function,
\begin{equation}
\label{eq:vertSup}
\begin{aligned}
h_\mathcal{B}(M^\top a^{\top}_{i,\mathcal{S}}) =  \underset{z}{\max} &~ \langle M^\top a^{\top}_{i,\mathcal{S}},z \rangle\\
                     \text{s.t.}          &~     A_\mathcal{B} z \leq b_\mathcal{B},
\end{aligned}
\end{equation}
where $A_{\mathcal{B}}=[\mathsf{I}_k~-\mathsf{I}_k]^\top$ and $b_{\mathcal{B}}=[\ol{v}^\top~\ul{v}^\top]^{\top}$. 
The above problem is an LP with a bounded feasible set. 
Thus, the optimal solution lies on the boundary of the hyper-rectangle towards which the normal $M^\top a^{\top}_{i,\mathcal{S}}$ points. Let us define, for all $i\in\Seq{1}{m}$, $
\hat{w}^i := \text{sign}(M^\top a^{\top}_{i,\mathcal{S}})\in\mathbb{R}^k$,
where the sign operator is applied entry-wise. 
(Notice that this vector simply indicates the orthant(s) that the vector $M^\top a^{\top}_{i,\mathcal{S}}$ points to.)
It then becomes clear that the vectors $w^i\in \mathbb{R}^{2k}$, as defined in \eqref{eq:NLP_form_w}, enable us to express the optimal solution of \eqref{eq:vertSup} in terms of a linear combination of the vertices of $\mathcal{B}$, i.e.,
\begin{align*}
h_{\mathcal{B}}\big( M^\top a^{\top}_{i,\mathcal{S}} \big)= \langle w^i, [\ol{v}^\top ~ \ul{v}^\top]^\top  \rangle, \forall i\in\Seq{1}{m}.
\end{align*}
Based on the above relation, the inequality~\eqref{eq:parBineq_1} simplifies to
\begin{align*}
\langle w^i, [\ol{v}^\top ~ \ul{v}^\top]^\top  \rangle \leq    b_{i,\mathcal{S}} - a_{i,\mathcal{S}} \xi, \forall i\in\Seq{1}{m},
\end{align*}
in which the vectors $\ul{v},\ol{v}\in\mathbb{R}^k$ are the decision variables. Intuitively, the above inequalities represents the linear constraints that the vertices of the hyper-rectangle $\mathcal{B}(-\ul{v},\ol{v})$ should satisfy in order to guarantee $\xi\in \mathcal{S}\ominus M \mathcal{B}(-\ul{v},\ol{v})$. 

Based on the chosen definition of volume for $\mathcal{B}(-\ul{v},\ol{v})$ in \eqref{eq:vol_def}, we intend to find a hyper-rectangle $\mathcal{B}(-\ul{v},\ol{v})$ that possesses the maximal volume. Unfortunately, regardless of the definition choice for the volume, the resulting objective function is non-convex and becomes unsuitable for optimization. 
Interestingly enough, one can simply use the logarithmic mapping for the volume definitions in \eqref{eq:vol_def} to obtain the objective functions suggested in \eqref{eq:cp_cost}, that are monotonic nonlinear concave functions. 
Then, it follows that a maximum hyper-rectangle $\mathcal{B}$ that contains the origin and satisfies $\xi\in \mathcal{S}\ominus M \mathcal{B}$ is the solution of the CP~\eqref{eq:NLP_form}.

\textbf{Theorem~\ref{theo:lp-nlp-Box_2}:} In the LP framework, 
%we look for the hyper-rectangle $\mathcal{B}\subset\mathbb{R}^k$ such that it has the maximal volume among all the hyper-rectangles contained in $\bar{\mathcal{S}}$ and $0\in\mathcal{B}$. Generally speaking, such a problem is called the ``{maximum containment problem}" that is an NP-hard problem. Instead 
we follow the procedure proposed in \cite{bemporad2004inner} with which one is able to cast the problem as a linear program. 
We first provide the proof for the LP relaxation of the problem~\eqref{eq:NLP_form} with $q=1$. 
Let us denote the maximum length of a line segment containing the origin, parallel to the $j$-th coordinate axis, and contained in $\bar{\mathcal{S}}$ by $r_j$. It follows from \cite[Proposition3]{bemporad2004inner} that one can use \eqref{eq:LP_rj} to find $r_j$, for all $j\in\Seq{1}{k}$. 
It is worth nothing that in the LP~\eqref{eq:LP_rj}, the constraints $z\leq 0$ and $z + \omega e_j \geq0$ are two extra regularity conditions that we placed on the line segment compared to \cite[Proposition3]{bemporad2004inner}. These conditions ensure that the origin lies inside this line segment.
Now, define the strictly positive vector $r\in\mathbb{R}^k$ by $r_j=\omega_j$ for all $j\in\Seq{1}{k}$. Then, it follows from \cite[Proposition2]{bemporad2004inner} that a maximum $r$-constrained inner hyper-rectangular $\mathcal{B}$ of $\bar{\mathcal{S}}$ that contains the origin is given by $\mathcal{B}(z^*,z^*+\lambda^*r)$ where $z^*$ and $\lambda^*$ are the optimal solutions of \eqref{eq:LP_form}.
Here, we also emphasize the fact that we have introduced the extra constraints $z\leq0$ and $z + \lambda r \geq0$ with respect to \cite[Proposition2]{bemporad2004inner}. By doing so, the LP~\eqref{eq:LP_form} is forced to find a hyper-rectangular $\mathcal{B}$ such that it contains the origin. Then, the claim for the LP case follows.

We now present a sketch of proof for the LP relaxation of the problem~\eqref{eq:NLP_form} with $q=2$. 
Observe that the polytope $\bar{\mathcal{S}}':=\big\{s\in\mathbb{R}^{2k}: W' s \leq B' \big\}$ is the inequality representation of the constraints in the CP~\eqref{eq:NLP_form}, where $W'$ and $B'$ are defined in Theorem~\ref{theo:lp-nlp-Box_2}. 
%Observe that the definition of volume~\eqref{eq:vol_def_2} intuitively suggests that we lift the polytope $\mathcal{S}$ to $\mathbb{R}^{2k}$ to promote a more symmetric hyper-rectangle $\mathcal{B}(l,u)$. 
We seek to find a hyper-rectangle that fits inside this \emph{lifted} polytope as follows. In the first step, we place a vertex of the hyper-rectangle at the origin.  
We then find the width of the line segment along each coordinate that is inside the lifted polytope and contains the origin using \eqref{eq:LP_rj_1}. 
In the second step, we use \eqref{eq:LP_form_1} to find a scaling factor $\lambda$ such that the $\lambda$-scaled hyper-rectangle constructed based on the first step fits inside the polytope $\bar{\mathcal{S}}'$. This concludes the proof.
%\end{proof}

\section{Numerical Examples}
\label{sec:exp}
In this section, we provide a numerical example to study the results presented in Section~\ref{sec:eb_RMPC}. For the numerical simulations, we use CVXOPT \cite{andersencvxopt} and (py)cddlib \cite{fukuda2001cddlib}. The system is an unstable batch reactor borrowed from \cite[Page~213]{rosenbrock1976computer}. We discretized the model using the zero-order-hold method with step-size~$0.05$, that is,   
\begin{equation*}
\begin{small}
x^+\!\!=\!\! \begin{pmatrix} 
 1.08 & -0.05 & 0.29 & -0.24\\
-0.03 & 0.81 & 0.00 & 0.03\\
 0.04 & 0.19 & 0.73 & 0.24\\
 0.00 & 0.19 & 0.05 & 0.91
\end{pmatrix}\!\! x +\!
\begin{pmatrix} 
 0.00 & -0.02\\
 0.26 & 0.00\\
 0.08 & -0.13\\
 0.08 & -0.00 
\end{pmatrix}\!\! u +w
\end{small}
\end{equation*}
where the state and input constraint sets are $\mathbb{X} =\{ x\in\mathbb{R}^4:~ \|x\|_{\infty} \leq 2 \}$ and $\mathbb{U} = \{ u\in\mathbb{R}^2:~ \| u \|_{\infty} \leq 2 \}$, respectively.
The disturbance set is defined as $\set{W}=\{ w\in\mathbb{R}^4:~\| w \|_{\infty} \leq 0.02\}$. 
The state and input target sets are $\mathbb{T}_x=\{ x\in\mathbb{R}^4:~ \| x\|_{\infty} \leq 0.5 \}$ and $\mathbb{T}_u= \{ u\in\mathbb{R}^2:~ \| u\|_{\infty} \leq 1.5 \}$, respectively.
The horizon length $N$ is set to $10$. The weight matrices in the cost function~\eqref{eq:cost_mpc} are $Q=2\times \mathsf{I}_4$ and $r=\mathsf{I}_2$. Finally, the terminal set is $\mathcal{X}_f=\{ x\in\mathbb{R}^4:~ \| x\|_{\infty} \leq 0.2 \}$.
%\begin{figure}[t]
%	\begin{center}
%		\centerline{\includegraphics[width=1\columnwidth]{fig1_inter_trig.pdf}}
%		\caption{Comparison of the CP and LP frameworks of Theorems~\ref{theo:lp-nlp-Box_1} \& \ref{theo:lp-nlp-Box_2} to construct the hyper-rectangles $\{\mathcal{E}_{j,k}\}_{j=1}^{N-1}$. The optimal state trajectory $\mathbf{X}^*_{k|k}=\{x^*_{k+i|k}\}_{i=0}^{N}$ is depicted by the (gray) solid/cross line. For several instance over the horizon $N=25$, the (colored) dash-dotted, dashed, dotted, solid rectangles represent the hyper-rectangles $\mathcal{E}_{j,k}$ using $\text{CP}_1$, $\text{LP}_1$, $\text{CP}_2$, and $\text{LP}_2$, respectively. For some $j\in\mathsf{N}_{[N-1]}$, the (colored) filled polytopes constructed around $x^*_{k+i|k}$ represent the principal polytopes $\bar{\mathcal{S}}_{j,k}$.}
%		\label{fig:comp_lp_nlp}
%	\end{center}
%\end{figure}

In what follows, we employ the triggering set construction approaches of Theorems~\ref{theo:lp-nlp-Box_1} \& \ref{theo:lp-nlp-Box_2} for $q=1$. 
Two types of disturbance realizations are considered: (1)~a uniform distribution with the bounded support~$\mathbb{W}$, and (2)~a worst case disturbance~$w_t = \text{argmax}_{w\in\mathbb{W}}~\xi_t^{\top}w$ at each time~$t$. 
In the case of uniform disturbance, we also manually applied an impulse-type disturbance to the closed-loop dynamics by resetting the second state~$\xi_{25}^2$ to $1.7$. This disturbance does not belong to the admissible disturbance set~$\mathbb{W}$.

\begin{figure}[t!]
	\centering
	\subfigure[Uniform case disturbance.]{\label{fig:CP1_uni}\includegraphics[width=.47\columnwidth]{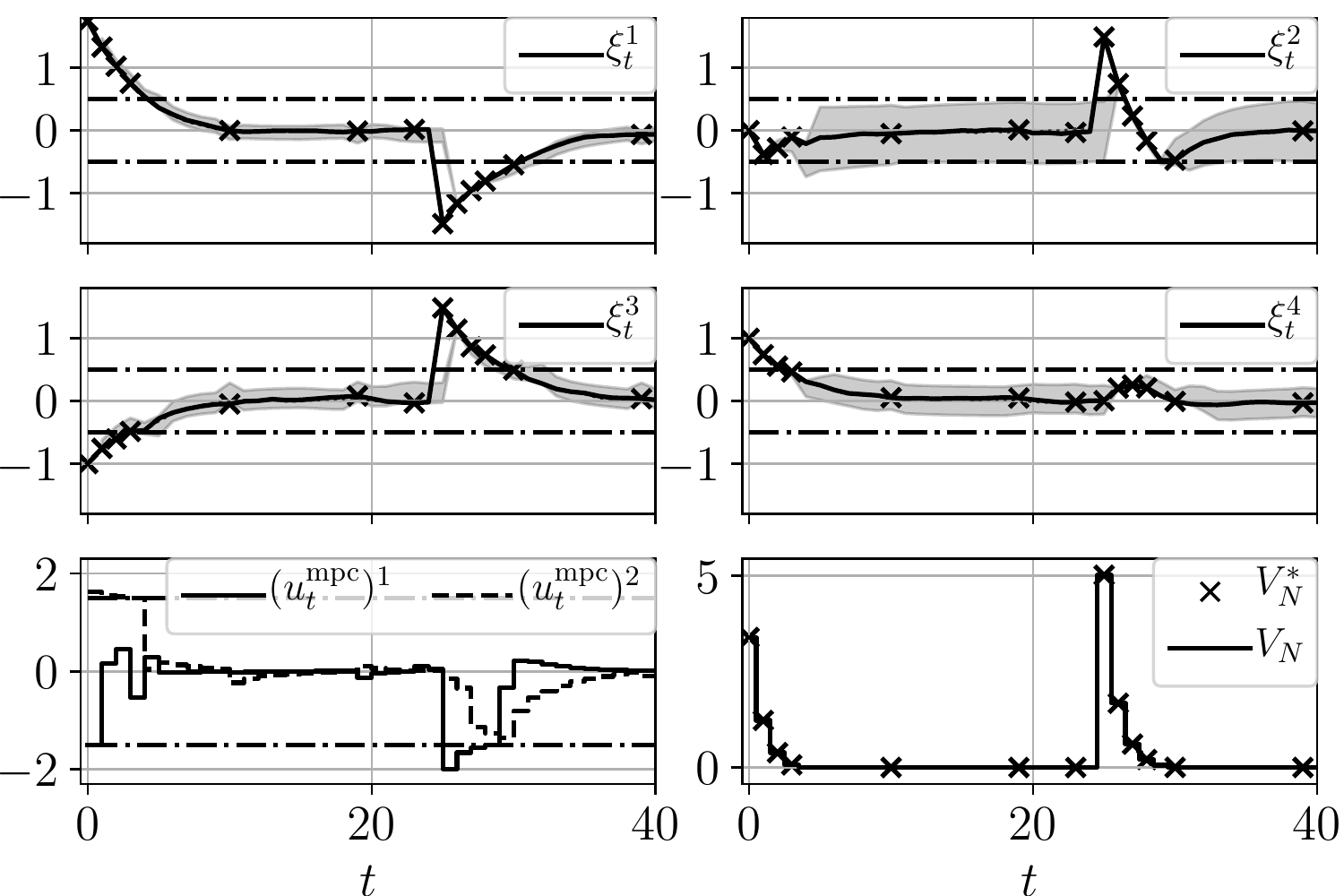}} 
	\subfigure[Worst case disturbance.]{\label{fig:CP1_worst}\includegraphics[width=.47\columnwidth]{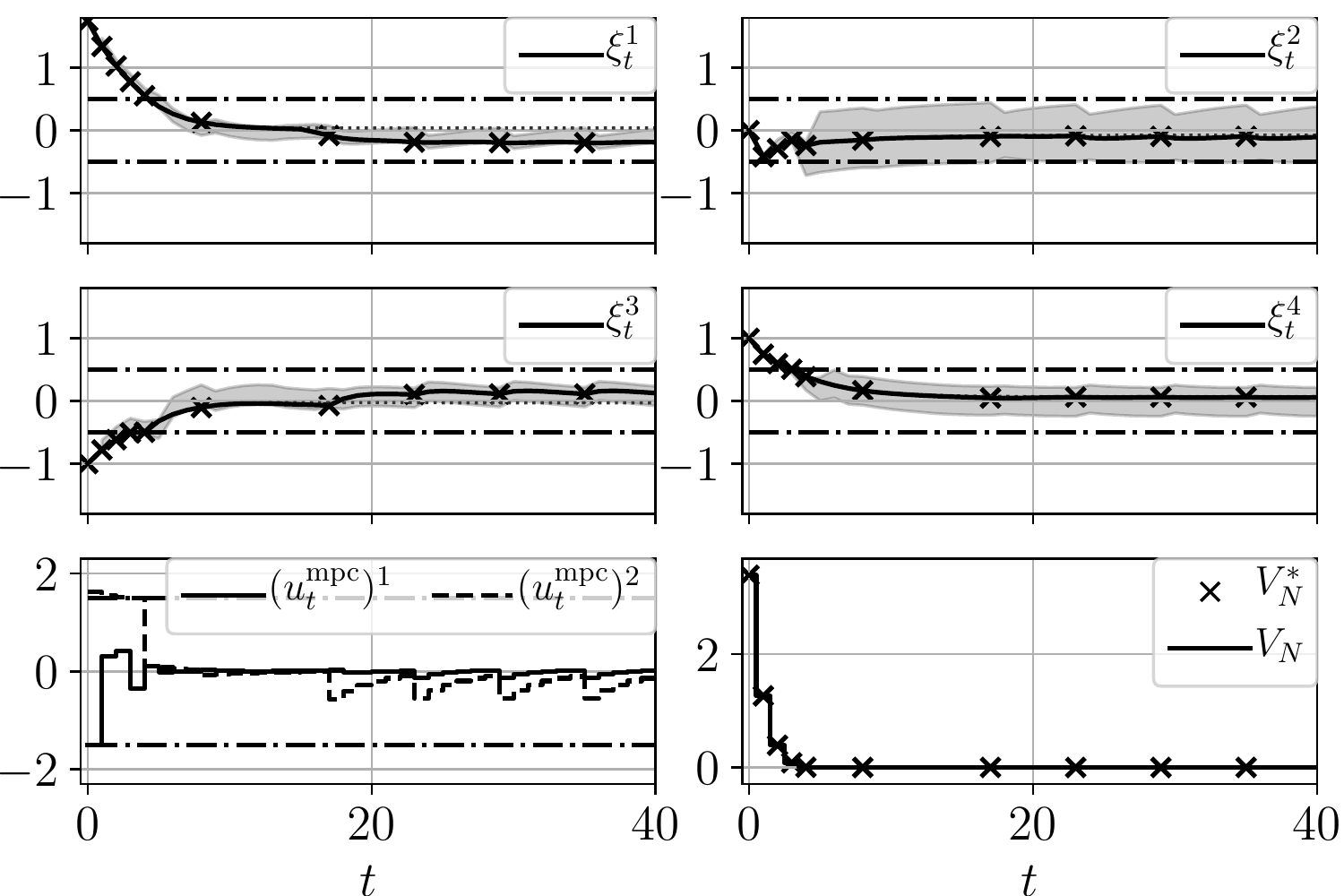}}
	\caption{Comparison of the event-based implementation using the construction approach $\text{CP}_1$ with the standard implementation. (Top four) The Solid lines are the evolution of states. The crosses are the states at triggering instances. The (gray) shaded areas are the projection of constructed hyper-rectangles~$\HypreC$ on the corresponding state's coordinate axis. (Bottom left) The lines are the input of the closed-loop system. (Bottom right) The crosses are the value function~$V^*_N$ at triggering instances. The solid line is the inter-event cost function~$V_N$ computed using Theorem~\ref{theo_4}.}
	\label{fig:CP1_uni_worst}
\end{figure}

\begin{figure}[t!]
	\centering
	\subfigure[Uniform case disturbance.]{\label{fig:LP1_uni}\includegraphics[width=.49\columnwidth]{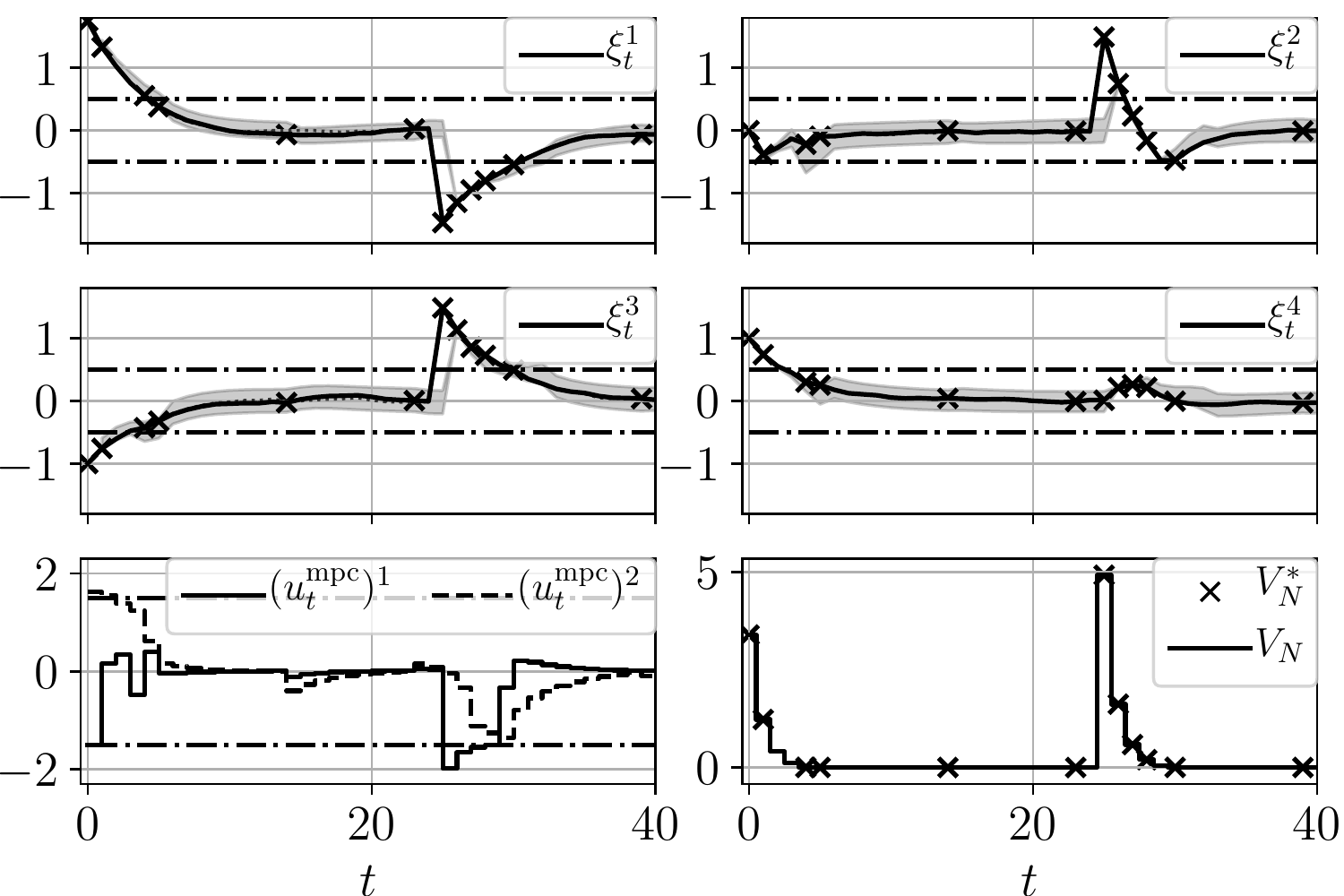}} 
	\subfigure[Worst case disturbance.]{\label{fig:LP1_worst}\includegraphics[width=.49\columnwidth]{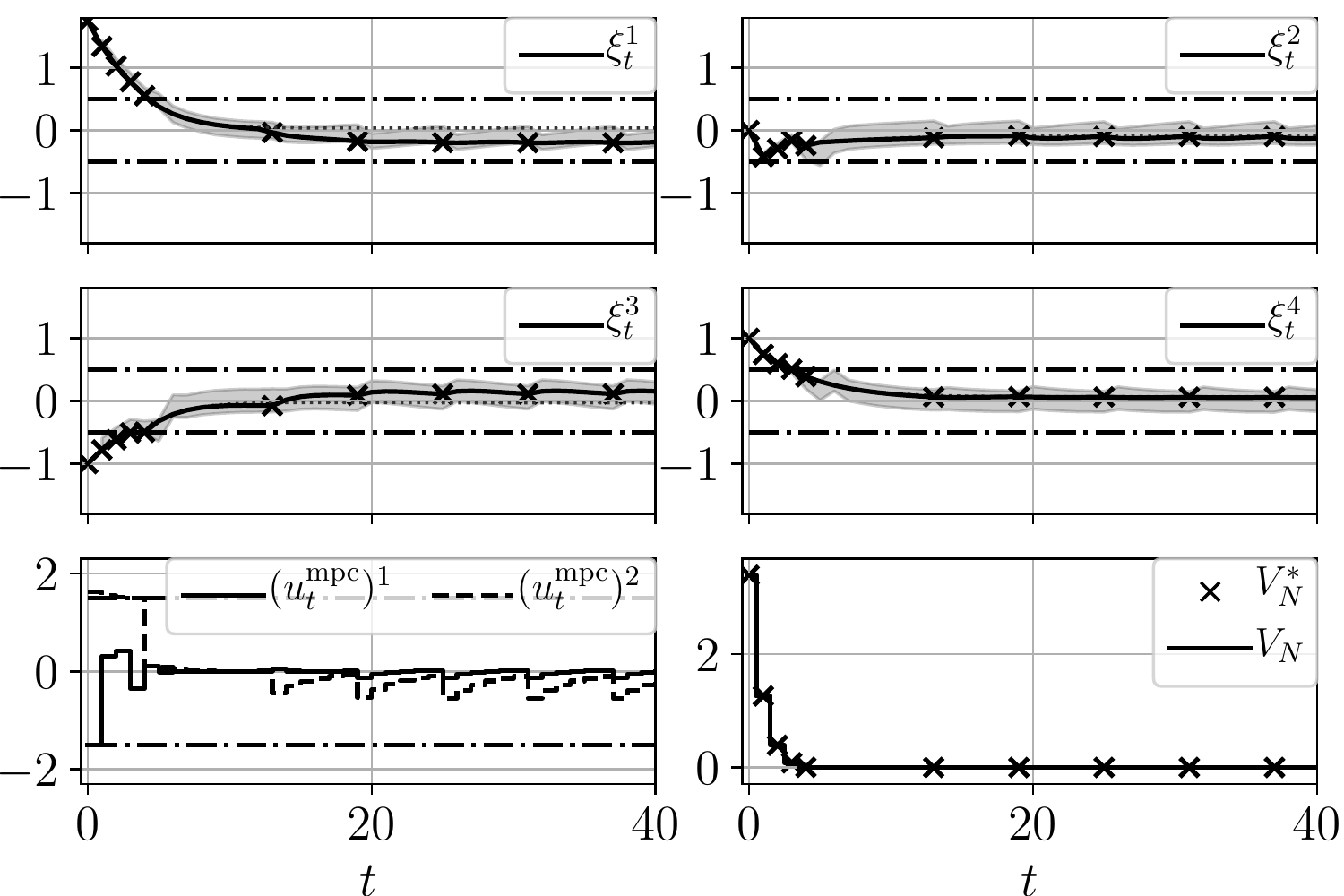}}
	\caption{Comparison of the event-based implementation using the construction approach $\text{LP}_1$ with the standard implementation. (Top four) The Solid lines are the evolution of states. The crosses are the states at triggering instances. The (gray) shaded areas are the projection of constructed hyper-rectangles~$\HypreC$ on the corresponding state's coordinate axis. (Bottom left) The lines are the input of the closed-loop system. (Bottom right) The crosses are the value function~$V^*_N$ at triggering instances. The solid line is the inter-event cost function~$V_N$ computed using Theorem~\ref{theo_4}.}
	\label{fig:LP1_uni_worst}
\end{figure}

Figures~\ref{fig:CP1_uni_worst} and \ref{fig:LP1_uni_worst} show the behavior of the event-based implementation of the MPC method. (Notice that the behavior of the standard MPC was almost identical, we did not include the results of the standard MPC for the sake clarity.) 

We begin with pointing out the shared properties of the approaches~$\text{CP}_1$ and $\text{LP}_1$. First of all, it is evident that the number of instances that the optimization problem~\eqref{MPC_tightened_target} is solved has effectively reduced in all considered cases compared to standard periodic implementation. 
Observe that the inputs and states of the closed-loop dynamics~\eqref{eq:closed_dyn} do not violate the constraint sets $\mathbb{X}$ and $\mathbb{U}$, respectively, in all considered cases. 
Moreover, the closed-loop states $\xi_t$ and the inputs $u_t$ converge to the target sets $\mathbb{T}_x$ and $\mathbb{T}_u$, respectively. Finally, both of the approaches~$\text{CP}_1$ and $\text{LP}_1$ can effectively recover from the impulse-type disturbance applied on time~$t=25$. 
We also note that the event-based implementations exhibit an almost limit-cyclic behavior inside the target set $\mathbb{T}_x$ in the worst case disturbance realizations. 

Let us now highlight the difference between the construction approaches~$\text{CP}_1$ and $\text{LP}_1$. 
As depicted in the top right plots of Figures~\ref{fig:CP1_uni} and \ref{fig:LP1_uni}, the construction method~$\text{LP}_1$ is more conservative in comparison with the construction method~$\text{CP}_1$. 
The width of the shaded areas represents the projection of the triggering sets~$\Hyprec$. 
In Figure~\ref{fig:CP1_uni}, one can also observe in the top right plot that the triggering intervals are tight with respect to the target sets, as well.

\section{Future Directions}
\label{sec:conc}
In this paper, an event-triggering approach was proposed to implement an RMPC method to constrained, perturbed LTI systems. The procedure to design the triggering mechanism is online, and is decoupled from the controller design. 
Specifically, we introduced two theoretical frameworks to construct the triggering mechanism as a volume maximization problem. 
%One framework is a general nonlinear convex program while the other framework is a linear program. 
%In particular, we proposed a non-standard definition of volume to address the limitations that occur in the case of using the standard definition of the volume in the process of designing the triggering mechanism. 
%For each choice of the volume definition, our numerical experiments showed that the theoretical frameworks provide a similar behavior at the price of the convex program framework being more computationally expensive compared to the linear program framework. 
%On the other hand, the linear and convex program frameworks based on the non-standard definition of the volume outperform the linear and convex program frameworks based on the standard definition of the volume. 
There are multiple directions that one can pursue to extend the results in this paper. 
First, it is interesting to investigate the possibility of extending the results of this paper to a nonlinear MPC case. 
In qualitative manner, we have observed that the choice of tightening gains~$\bm{K}$ directly impacts the constructed triggering sets~$\bm{\mathcal{E}}$. Hence, another possible direction is to explore the possibility of characterizing this unknown dependency in a more quantitative manner. Lastly, the triggering approach proposed in this paper is online (and in fact state-dependent). It is thus valuable to investigate whether it is possible to make the triggering set design offline.

%\begin{figure}[t]
%	\begin{center}
%		\centerline{\includegraphics[width=.7\columnwidth]{fig1_inter_trig.pdf}}
%		\caption{Comparison of the CP and LP frameworks of Theorems~\ref{theo:lp-nlp-Box_1} \& \ref{theo:lp-nlp-Box_2} to construct the hyper-rectangles $\{\mathcal{E}_{j,k}\}_{j=1}^{N-1}$. The optimal state trajectory $\mathbf{X}^*_{k|k}=\{x^*_{k+i|k}\}_{i=0}^{N}$ is depicted by the (gray) solid/cross line. For several instant over the horizon $N=25$, the (colored) dash-dotted, dashed, dotted, solid rectangles represent the hyper-rectangles $\mathcal{E}_{j,k}$ using $\text{CP}_1$, $\text{LP}_1$, $\text{CP}_2$, and $\text{LP}_2$, respectively. For some $j\in\mathsf{N}_{[N-1]}$, the (colored) filled polytopes constructed around $x^*_{k+i|k}$ represent the principal polytopes $\bar{\mathcal{S}}_{j,k}$.}
%		\label{fig:comp_lp_nlp}
%	\end{center}
%\end{figure}

%\newpage
%===============================================================================
	\bibliographystyle{siam}	
	\bibliography{./mybref_ETRMPC}

\end{document}